# MARTINGALE APPROACH TO STOCHASTIC DIFFERENTIAL GAMES OF CONTROL AND STOPPING[1]

By Ioannis Karatzas and Ingrid-Mona Zamfirescu

*Columbia University and Baruch College, CUNY*

We develop a martingale approach for studying continuous-time stochastic differential games of control and stopping, in a non-Markovian framework and with the control affecting only the drift term of the state-process. Under appropriate conditions, we show that the game has a value and construct a *saddle pair* of optimal control and stopping strategies. Crucial in this construction is a characterization of saddle pairs in terms of pathwise and martingale properties of suitable quantities.

**1. Introduction and synopsis.** We develop a theory for zero-sum stochastic differential games with two players, a "controller" and a "stopper." The state $X(\cdot)$ in these games evolves in Euclidean space according to a stochastic functional/differential equation driven by a Wiener process; via his choice of instantaneous, nonanticipative control $u(t)$, the controller can affect the local drift of this state process $X(\cdot)$ at time $t$, though not its local variance.

The stopper decides the duration of the game, in the form of a stopping rule $\tau$ for the process $X(\cdot)$. At the terminal time $\tau$ the stopper receives from the controller a "reward" $\int_0^\tau h(t, X, u_t) \, dt + g(X(\tau))$ consisting of two parts: The integral up to time $\tau$ of a time-dependent *running reward* $h$, which also depends on the past and present states $X(s), 0 \le s \le t$, and on the present value $u_t$ of the control; and the value at the terminal state $X(\tau)$ of a continuous *terminal reward* function $g$ ("reward" always refers to the stopper).

Under appropriate conditions on the local drift and local variance of the state process, and on the running and terminal cost functions $h$ and $g$,

Received December 2006; revised June 2007.
[1]Supported in part by the NSF Grant DMS-06-01774.

we establish the existence of a value for the resulting stochastic game of control and stopping, as well as regularity and martingale-type properties of the temporal evolution for the resulting *value process*. We also construct optimal strategies for the two players, in the form of a *saddle point* $(u^*, \tau_*)$, to wit: the strategy $u^*(\cdot)$ is the controller's best response to the stopper's use of the stopping rule $\tau_*$, in the sense of minimizing total expected cost; and the stopping rule $\tau_*$ is the stopper's best response to the controller's use of the control strategy $u^*(\cdot)$, in the sense of maximizing total expected reward.

The approach of the paper is direct and probabilistic. It draws on the Dubins–Savage (1965) theory, and builds on the martingale methodologies developed for the optimal stopping problem and for the problem of optimal stochastic control over the last three decades; see, for instance, Neveu (1975), El Karoui (1981), Beneš (1970, 1971), Rishel (1970), Duncan and Varaiya (1971), Davis and Varaiya (1973), Davis (1973, 1979) and Elliott (1977, 1982). It proceeds in terms of a characterization of saddle points via martingale-type properties of suitable quantities, which involve the value process of the game.

An advantage of the approach is that it imposes no Markovian assumptions on the dynamics of the state-process; it allows the local drift and variance of the state-process, as well as the running cost, to depend at any given time $t$ on past-and-present states $X(s), 0 \leq s \leq t$, in a fairly general, measurable manner. (The boundedness and continuity assumptions can most likely be relaxed.)

The main drawback of this approach is that it imposes a severe nondegeneracy condition on the local variance of the state-process, and does not allow this local variance to be influenced by the controller. We hope that subsequent work will be able to provide a more general theory for such stochastic games, possibly also for their nonzero-sum counterparts, without such restrictive assumptions—at least in the Markovian framework of, say, Fleming and Soner (2006), El Karoui, Nguyen and Jeanblanc-Picqué (1987), Bensoussan and Lions (1982) or Bismut (1973, 1978). It would also be of considerable interest to provide a theory for control of "bounded variation" type (admixture of absolutely continuous, as in this paper, with pure jump and singular, terms).

**Extant work:** A game between a controller and a stopper, in discrete time and with Polish (complete separable metric) state-space, was studied by Maitra and Sudderth (1996b); under appropriate conditions, these authors obtained the existence of a value for the game and provided a transfinite induction algorithm for its computation.

In Karatzas and Sudderth (2001) a similar game was studied for a linear diffusion process, with the unit interval as its state-space and absorption at the endpoints. The one-dimensional nature of the setup allowed an explicit computation of the value and of a saddle pair of strategies, based on scale-function considerations and under a nondegeneracy condition on the variance of the diffusion. Karatzas and Sudderth (2007) studied recently nonzero-sum versions of these linear diffusion games, where one seeks and constructs Nash equilibria, rather than saddle pairs, of strategies. Always in a Markovian, one-dimensional framework, Weerasinghe (2006) was able to solve in a similar, explicit manner, a stochastic game with variance that *is* allowed to degenerate; while Bayraktar and Young (2007) established a very interesting convex-duality connection, between a stochastic game of control and stopping and a probability-of-ruin-minimization problem.

Along a parallel tack, stochastic games of stopping have been treated via the theory of *Backwards Stochastic Differential Equations* starting with Cvitanić and Karatzas (1996), and continuing with Hamadène and Lepeltier (1995, 2000) and Hamadène (2006) for games of mixed control/stopping.

The methods used in the present paper are entirely different from those in all these works.

• The coöperative version of the game has received far greater attention. In the standard model of stochastic control, treated, for instance, in the classic monograph Fleming and Soner (1992), the controller may influence the state dynamics but must operate over a prescribed time-horizon. If the controller is also allowed to choose a quitting time adaptively, at the expense of incurring a termination cost, one has a problem of *control with discretionary stopping* [or "leavable" control problem, in the terminology of Dubins and Savage (1976)]. General existence/characterization results for such problems were obtained by Dubins and Savage (1976) and by Maitra and Sudderth (1996a) under the rubric of "leavable" stochastic control; by Krylov (1980), El Karoui (1981), Bensoussan and Lions (1982), Haussmann and Lepeltier (1990), Maitra and Sudderth (1996a), Morimoto (2003), Ceci and Basan (2004); and by Karatzas and Zamfirescu (2006) in the present framework. There are also several explicitly solvable problems in this vein: see Beneš (1992), Davis and Zervos (1994), Karatzas and Sudderth (1999), Karatzas et al. (2000), Karatzas and Wang (2000, 2001), Kamizono and Morimoto (2002), Karatzas and Ocone (2002), Ocone and Weerasinghe (2006).

Such problems arise, for instance, in target-tracking models, where one has to stay close to a target by spending fuel, declare when one has arrived "sufficiently close," then decide whether to engage the target or not. Combined stochastic control/optimal stopping problems also arise in mathematical finance, namely, in the context of computing the upper-hedging prices of *American contingent claims under constraints*; these computations lead to stochastic control of the absolutely continuous or the singular type [e.g. Karatzas and Kou (1998), Karatzas and Wang (2000)].

The computation of the lower-hedging prices for American contingent claims under constraints leads to stochastic games of control and stopping; see Karatzas and Kou (1998) for details.

**Synopsis:** We set up in the next section the model for a controlled stochastic functional/differential equation driven by a Wiener process, that will be used throughout the paper; this setting is identical to that of Elliott (1982) and of our earlier paper Karatzas and Zamfirescu (2006). Within this model, we formulate in Section 3 the stochastic game of control and stopping that will be the focus of our study. Section 4 reviews in the present context the classical results for *optimal stopping* on the one hand, and for *optimal stochastic control* on the other, when these problems are viewed separately.

Section 5 establishes the existence of a value for the stochastic game, and studies the regularity and some simple martingale-like properties of the resulting value process evolving through time. This study continues in earnest in Section 6 and culminates with Theorem 6.3.

Section 7 then builds on these results, to provide necessary and sufficient conditions for a pair $(u, \tau)$ consisting of a control strategy and a stopping rule, to be a saddle point for the stochastic game. These conditions are couched in terms of martingale-like properties for suitable quantities, which involve the value process and the cumulative running reward. A similar characterization is provided in Section 8 for the optimality of a given control strategy $u(\cdot)$.

With the help of the predictable representation property of the Brownian filtration under equivalent changes of probability measure, and of the Doob–Meyer decomposition for sufficiently regular submartingales, this characterization leads—in Section 9, and under appropriate conditions—to a specific control strategy $u^*(\cdot)$ as a candidate for optimality. These same martingale-type conditions suggest $\tau_*$, the first time the value process $V(\cdot)$ of the game agrees with the terminal reward $g(X(\cdot))$, as a candidate for optimal stopping rule. Finally, it is shown that the pair $(u^*, \tau_*)$ is indeed a saddle point of the stochastic game, and that $V(\cdot \wedge \tau_*)$ has continuous paths.

NOTATION.    The paper is quite heavy with notation, so here is a partial list for ease of reference:

$X(t), W^u(t)$: Equations (1), (6) and equation (4), respectively.

$\Lambda^u(\mathfrak{t}), \Lambda^u(\mathfrak{t}, \tau)$: Exponential likelihood ratios (martingales); equations (3) and (25).

$Y^u(\mathfrak{t}, \tau), Y^u(\tau)$: Total (i.e., terminal, plus running) cost/reward on the interval $[[\mathfrak{t}, \tau]]$: equations (8), (23).

$\overline{V}, \underline{V}$ and $\overline{V}(\mathfrak{t}), \underline{V}(\mathfrak{t})$: Upper and lower values of the game; equations (9), (11), (12).

$J(\mathfrak{t}, \tau)$: Minimal conditional expected total cost on the interval $[[\mathfrak{t}, \tau]]$; equation (14).

$Z^u(\mathfrak{t})$: Maximal conditional expected reward under control $u(\cdot)$, from time $\mathfrak{t}$ onward; equation (19).
$Q^u(\mathfrak{t})$: Cumulative maximal conditional expected reward under control $u(\cdot)$, at time $\mathfrak{t}$; equation (20).
$R^u(\mathfrak{t})$: Cumulative value of game under control $u(\cdot)$, at $\mathfrak{t}$; equation (36).
$\tau_\mathfrak{t}^u(\varepsilon), \tau_\mathfrak{t}^u$: Stopping rules; equation (22).
$\varrho_\mathfrak{t}(\varepsilon), \varrho_\mathfrak{t}$: Stopping rules; equation (33).
$H(t, \omega, a, p)$: Hamiltonian function; equation (72).
*Saddle Point*: Inequalities (10).
*Thrifty Control Strategy*: Requirement (66).

**2. The model.** Consider the space $\Omega = C([0,T]; \mathbb{R}^n)$ of continuous functions $\omega : [0,T] \to \mathbb{R}^n$, defined on a given bounded interval $[0,T]$ and taking values in some Euclidean space $\mathbb{R}^n$. The coördinate mapping process will be denoted by $W(t, \omega) = \omega(t), 0 \le t \le T$, and $\mathcal{F}_t^W = \sigma(W(s); 0 \le s \le t)$, $0 \le t \le T$, will stand for the natural filtration generated by this process $W$. The measurable space $(\Omega, \mathcal{F}_T^W)$ will be endowed with Wiener measure $\mathbb{P}$, under which $W$ becomes a standard $n$-dimensional Brownian motion. We shall denote by $\mathbb{F} = \{\mathcal{F}_t\}_{0 \le t \le T}$ the $\mathbb{P}$-augmentation of this natural filtration, and use the notation $\|\omega\|_t^* := \max_{0 \le s \le t} |\omega(s)|$, $\omega \in \Omega, 0 \le t \le T$.

The $\sigma$-algebra of predictable subsets of the product space $[0,T] \times \Omega$ will be denoted by $\mathcal{P}$, and $\mathcal{S}$ will stand for the collection of *stopping rules* of the filtration $\mathbb{F}$. These are measurable mappings $\tau : \Omega \to [0,T]$ with the property

$$\{\tau \le t\} \in \mathcal{F}_t \qquad \forall \, 0 \le t \le T.$$

Given any two stopping rules $\rho$ and $\nu$ with $\rho \le \nu$, we shall denote by $\mathcal{S}_{\rho, \nu}$ the collection of all stopping rules $\tau \in \mathcal{S}$ with $\rho \le \tau \le \nu$.

Consider now a *predictable* (i.e., $\mathcal{P}$-measurable) $\sigma : [0,T] \times \Omega \to \mathbb{L}(\mathbb{R}^n; \mathbb{R}^n)$ with values in the space $\mathbb{L}(\mathbb{R}^n; \mathbb{R}^n)$ of $(n \times n)$ matrices, and suppose that $\sigma(t, \omega)$ is nonsingular for every $(t, \omega) \in [0,T] \times \Omega$ and that there exists some real constant $K > 0$ for which

$$\|\sigma^{-1}(t, \omega)\| \le K \quad \text{and} \quad |\sigma_{ij}(t, \omega) - \sigma_{ij}(t, \widetilde{\omega})| \le K \|\omega - \widetilde{\omega}\|_t^* \qquad \forall \, 1 \le i, j \le n,$$

hold for every $\omega \in \Omega, \widetilde{\omega} \in \Omega$ and every $t \in [0,T]$. Then for any initial condition $x \in \mathbb{R}^n$, there is a pathwise unique, strong solution $X(\cdot)$ of the stochastic equation

$$(1) \qquad X(t) = x + \int_0^t \sigma(s, X) \, dW(s), \qquad 0 \le t \le T;$$

see Theorem 14.6 in Elliott (1982). In particular, the augmentation of the natural filtration generated by $X(\cdot)$ coincides with the filtration $\mathbb{F}$ itself.

Now let us introduce an element of control in this picture. We shall denote by $\mathfrak{U}$ the class of *admissible control strategies* $u : [0,T] \times \Omega \to A$. These are predictable processes with values in some given separable metric space $A$. We shall assume that $A$ is a countable union of nonempty, compact subsets, and is endowed with the $\sigma$-algebra $\mathcal{A}$ of its Borel subsets.

We shall consider also a $\mathcal{P} \otimes \mathcal{A}$-measurable function $f:([0,T] \times \Omega) \times A \to \mathbb{R}^n$ with the following properties:

- for each $a \in A$, the mapping $(t, \omega) \mapsto f(t, \omega, a)$ is predictable; and
- there exists a real constant $K > 0$ such that

$$(2) \qquad |f(t,\omega,a)| \leq K(1+\|\omega\|_t^*) \qquad \forall\, 0 \leq t \leq T, \omega \in \Omega, a \in A.$$

For any given admissible control strategy $u(\cdot) \in \mathfrak{U}$, the exponential process

$$(3) \qquad \begin{aligned} \Lambda^u(t) := \exp\bigg\{ & \int_0^t \langle \sigma^{-1}(s,X) f(s,X,u_s), dW(s) \rangle \\ & - \tfrac{1}{2} \int_0^t |\sigma^{-1}(s,X) f(s,X,u_s)|^2\, ds \bigg\} \end{aligned}$$

$0 \leq t \leq T$, is a martingale under all these assumptions; namely, $\mathbb{E}(\Lambda^u(T)) = 1$ [see Beneš (1971), as well as Karatzas and Shreve (1991), pages 191 and 200 for this result]. Then the Girsanov theorem (ibid., Section 3.5) guarantees that the process

$$(4) \qquad W^u(t) := W(t) - \int_0^t \sigma^{-1}(s,X) f(s,X,u_s)\, ds, \qquad 0 \leq t \leq T$$

is a Brownian motion with respect to the filtration $\mathbb{F}$, under the new probability measure

$$(5) \qquad \mathbb{P}^u(B) := \mathbb{E}[\Lambda^u(T) \cdot 1_B], \qquad B \in \mathcal{F}_T,$$

which is equivalent to $\mathbb{P}$. It is now clear from the equations (1) and (4) that

$$(6) \quad X(t) = x + \int_0^t f(s,X,u_s)\, ds + \int_0^t \sigma(s,X)\, dW^u(s), \qquad 0 \leq t \leq T,$$

holds almost surely. This will be our model for a controlled stochastic functional/differential equation, with the control appearing only in the drift (bounded variation) term.

**3. The stochastic game of control and stopping.** In order to specify the objective of our stochastic game of control and stopping, let us consider two bounded, measurable functions $h:[0,T] \times \Omega \times A \to \mathbb{R}$ and $g:\mathbb{R}^n \to \mathbb{R}$. We shall assume that the *running reward* function $h$ satisfies the measurability conditions imposed on the drift-function $f$ above, except of course that (2) is now strengthened to the boundedness requirement

$$(7) \qquad |h(t,\omega,a)| \leq K \qquad \forall\, 0 \leq t \leq T,\ \omega \in \Omega,\ a \in A.$$

To simplify the analysis, we shall assume that the *terminal reward* function $g$ is continuous.

We shall study a stochastic game of control and stopping with two players: The *controller*, who chooses an admissible control strategy $u(\cdot)$ in $\mathfrak{U}$; and the *stopper*, who decides the duration of the game by his choice of stopping rule $\tau \in \mathcal{S}$. When the stopper declares the game to be over, he receives from the controller the amount $Y^u(\tau) \equiv Y^u(0,\tau)$, where

$$(8) \qquad Y^u(\mathfrak{t},\tau) := g(X(\tau)) + \int_{\mathfrak{t}}^{\tau} h(s, X, u_s)\, ds \qquad \text{for } \tau \in \mathcal{S}_{\mathfrak{t},T},\ \mathfrak{t} \in \mathcal{S}.$$

It is thus in the best interest of the controller (resp., the stopper) to make the amount $Y^u(\tau)$ as small (resp., as large) as possible, at least on the average. We are thus led to a *stochastic game*, with

$$(9) \qquad \overline{V} := \inf_{u \in \mathfrak{U}} \sup_{\tau \in \mathcal{S}} \mathbb{E}^u(Y^u(\tau)), \qquad \underline{V} := \sup_{\tau \in \mathcal{S}} \inf_{u \in \mathfrak{U}} \mathbb{E}^u(Y^u(\tau))$$

as its upper- and lower-values, respectively; clearly, $\underline{V} \leq \overline{V}$.

We shall say that *the game has a value*, if its upper- and lower-values coincide, that is, $\underline{V} = \overline{V}$; in that case we shall denote this common value simply by $V$.

A pair $(u^*, \tau_*) \in \mathfrak{U} \times \mathcal{S}$ will be called *saddle point* of the game, if

$$(10) \qquad \mathbb{E}^{u^*}(Y^{u^*}(\tau)) \leq \mathbb{E}^{u^*}(Y^{u^*}(\tau_*)) \leq \mathbb{E}^u(Y^u(\tau_*))$$

holds for every $u(\cdot) \in \mathfrak{U}$ and $\tau \in \mathcal{S}$. In words, the strategy $u^*(\cdot)$ is the controller's best response to the stopper's use of the rule $\tau_*$; and the rule $\tau_*$ is the stopper's best response to the controller's use of the strategy $u^*(\cdot)$.

If such a saddle-point pair $(u^*, \tau^*)$ exists, then it is quite clear that the game has a value. We shall characterize the saddle property in terms of simple, pathwise and martingale properties of certain crucial quantities; see Theorem 7.1. Then, in Sections 8 and 9, we shall use this characterization in an effort to show that a saddle point indeed exists and to identify its components.

In this effort we shall need to consider, a little more generally than in (9), the *upper-value-process*

$$(11) \qquad \overline{V}(\mathfrak{t}) := \operatorname*{essinf}_{u \in \mathfrak{U}} \operatorname*{esssup}_{\tau \in \mathcal{S}_{\mathfrak{t},T}} \mathbb{E}^u(Y^u(\mathfrak{t},\tau)|\mathcal{F}_{\mathfrak{t}})$$

and the *lower-value-process*

$$(12) \qquad \underline{V}(\mathfrak{t}) := \operatorname*{esssup}_{\tau \in \mathcal{S}_{\mathfrak{t},T}} \operatorname*{essinf}_{u \in \mathfrak{U}} \mathbb{E}^u(Y^u(\mathfrak{t},\tau)|\mathcal{F}_{\mathfrak{t}})$$

of the game, for each $\mathfrak{t} \in \mathcal{S}$. Clearly, $\overline{V}(0) = \overline{V}$, $\underline{V}(0) = \underline{V}$, as well as

$$(13) \qquad g(X(\mathfrak{t})) \leq \underline{V}(\mathfrak{t}) \leq \overline{V}(\mathfrak{t}) \qquad \forall\, \mathfrak{t} \in \mathcal{S}.$$

We shall see in Theorem 5.1 that this last inequality holds, in fact, as an equality: the game has a value at all times.

**4. Optimal control and stopping problems, viewed separately.** Given any stopping rule $\mathfrak{t} \in \mathcal{S}$, we introduce the minimal conditional expected cost

$$
(14) \qquad J(\mathfrak{t}, \tau) := \operatorname*{essinf}_{v \in \mathfrak{U}} \mathbb{E}^v(Y^v(\mathfrak{t}, \tau) | \mathcal{F}_{\mathfrak{t}}),
$$

that can be achieved by the controller over the stochastic interval

$$
(15) \qquad [[\mathfrak{t}, \tau]] := \{(s, \omega) \in [0, T] \times \Omega : \mathfrak{t}(\omega) \le s \le \tau(\omega)\},
$$

for each stopping rule $\tau \in \mathcal{S}_{\mathfrak{t},T}$. With the notation (14), the lower value (12) of the game becomes

$$
(16) \qquad \underline{V}(\mathfrak{t}) = \operatorname*{esssup}_{\tau \in \mathcal{S}_{\mathfrak{t},T}} J(\mathfrak{t}, \tau) \ge J(\mathfrak{t}, \mathfrak{t}) = g(X(\mathfrak{t})) \qquad \text{a.s.}
$$

By analogy with the classical martingale approach to stochastic control [developed by Rishel (1970), Duncan and Varaiya (1971), Davis and Varaiya (1973), Davis (1973) and outlined in Davis (1979), El Karoui (1981)], for any given admissible control strategy $u(\cdot) \in \mathfrak{U}$ and any stopping rules $\mathfrak{t}, \nu, \tau$ with $0 \le \mathfrak{t} \le \nu \le \tau \le T$, we have the $\mathbb{P}^u$-submartingale property

$$
(17) \qquad \mathbb{E}^u(\Psi^u(\nu, \tau) | \mathcal{F}_{\mathfrak{t}}) \ge \Psi^u(\mathfrak{t}, \tau)
$$

$$
\text{for } \Psi^u(\mathfrak{t}, \tau) := J(\mathfrak{t}, \tau) + \int_0^{\mathfrak{t}} h(s, X, u_s)\, ds,
$$

or equivalently,

$$
(18) \qquad \mathbb{E}^u \left[ J(\nu, \tau) + \int_{\mathfrak{t}}^{\nu} h(s, X, u_s)\, ds \Big| \mathcal{F}_{\mathfrak{t}} \right] \ge J(\mathfrak{t}, \tau) \qquad \text{a.s.}
$$

A very readable account of this theory appears in Chapter 16, pages 222–241 of Elliott (1982).

4.1. *A family of optimal stopping problems.* For each admissible control strategy $u(\cdot) \in \mathfrak{U}$, we define the maximal conditional expected reward

$$
(19) \qquad Z^u(\mathfrak{t}) := \operatorname*{esssup}_{\tau \in \mathcal{S}_{\mathfrak{t},T}} \mathbb{E}^u(Y^u(\mathfrak{t}, \tau) | \mathcal{F}_{\mathfrak{t}}), \qquad \mathfrak{t} \in \mathcal{S},
$$

that can be achieved by the stopper from time $\mathfrak{t}$ onward, as well as the "cumulative" quantity

$$
(20) \qquad Q^u(\mathfrak{t}) := Z^u(\mathfrak{t}) + \int_0^{\mathfrak{t}} h(s, X, u_s)\, ds = \operatorname*{esssup}_{\tau \in \mathcal{S}_{\mathfrak{t},T}} \mathbb{E}^u(Y^u(\tau) | \mathcal{F}_{\mathfrak{t}});
$$

in particular,

$$
(21) \qquad Z^u(\mathfrak{t}) \ge Y^u(\mathfrak{t}, \mathfrak{t}) = g(X(\mathfrak{t})), \qquad \overline{V}(\mathfrak{t}) = \operatorname*{essinf}_{u \in \mathfrak{U}} Z^u(\mathfrak{t}).
$$

Let us introduce the stopping rules

$$(22) \quad \tau_{\mathfrak{t}}^u(\varepsilon) := \inf\{s \in [\mathfrak{t}, T] : g(X(s)) \geq Z^u(s) - \varepsilon\}, \qquad \tau_{\mathfrak{t}}^u := \tau_{\mathfrak{t}}^u(0)$$

for each $\mathfrak{t} \in \mathcal{S}, 0 \leq \varepsilon < 1$. Then $\tau_{\mathfrak{t}}^u(\varepsilon) \leq \tau_{\mathfrak{t}}^u$.

From the classical martingale approach to the theory of optimal stopping [e.g., El Karoui (1981) or Karatzas and Shreve (1998), Appendix D], the following results are well known.

PROPOSITION 4.1. *The process $Q^u(\cdot) \equiv \{Q^u(t), 0 \leq t \leq T\}$ is a $\mathbb{P}^u$-supermartingale with paths that are RCLL (Right-Continuous, with Limits from the Left); it dominates the continuous process $Y^u(\cdot)$ given as*

$$(23) \quad Y^u(t) \equiv Y^u(0, t) = g(X(t)) + \int_0^t h(s, X, u_s)\, ds, \qquad 0 \leq t \leq T;$$

*and $Q^u(\cdot)$ is the smallest RCLL supermartingale which dominates $Y^u(\cdot)$.*

In other words, $Q^u(\cdot)$ of (20) is the *Snell Envelope* of the process $Y^u(\cdot)$.

PROPOSITION 4.2. *For any stopping rules $\mathfrak{t}, \nu, \theta$ with $\mathfrak{t} \leq \nu \leq \theta \leq \tau_{\mathfrak{t}}^u$, we have the martingale property $\mathbb{E}^u[Q^u(\theta)|\mathcal{F}_\nu] = Q^u(\nu)$ a.s.; in particular, $Q^u(\cdot \wedge \tau_0^u)$ is a $\mathbb{P}^u$-martingale. Furthermore, $Z^u(\mathfrak{t}) = \mathbb{E}^u[Y^u(\mathfrak{t}, \tau_{\mathfrak{t}}^u)|\mathcal{F}_\mathfrak{t}]$ holds a.s.*

4.2. *A preparatory lemma.* For the proof of several results in this work, we shall need the following observation; we list it separately, for ease of reference.

LEMMA 4.3. *Suppose that $\mathfrak{t}, \theta$ are stopping rules with $0 \leq \mathfrak{t} \leq \theta \leq T$, and that $u(\cdot), v(\cdot)$ are admissible control strategies in $\mathfrak{U}$.*

(i) *Assume that $u(\cdot) = v(\cdot)$ holds a.e. on the stochastic interval $[[\mathfrak{t}, \theta]]$, in the notation of (15). Then, for any bounded and $\mathcal{F}_\theta$-measurable random variable random variable $\Xi$, we have*

$$(24) \qquad \mathbb{E}^v[\Xi|\mathcal{F}_\mathfrak{t}] = \mathbb{E}^u[\Xi|\mathcal{F}_\mathfrak{t}] \qquad \text{a.s.}$$

*In particular, with $\mathfrak{t} = 0$ this gives $\mathbb{E}^v[\Xi] = \mathbb{E}^u[\Xi]$.*

(ii) *More generally, assume that $u(\cdot) = v(\cdot)$ holds a.e. on $\{(u, \omega) : \mathfrak{t}(\omega) \leq u \leq \theta(\omega), \omega \in A\}$ for some $A \in \mathcal{F}_\mathfrak{t}$. Then (24) holds a.e. on the event $A$.*

The reasoning is simple: with the notation $\Lambda^u(\mathfrak{t}, \theta) := \Lambda^u(\theta)/\Lambda^u(\mathfrak{t})$ from (3), and using the martingale property of $\Lambda^u(\cdot)$ under $\mathbb{P}^u$, we have $\mathbb{E}^u[\Lambda^u(\mathfrak{t}, \theta)|\mathcal{F}_\mathfrak{t}] =$ 1 a.s. In conjunction with the Bayes rule for conditional expectations under equivalent probability measures, this gives

$$\mathbb{E}^u[\Xi|\mathcal{F}_{\mathfrak{t}}] = \frac{\Lambda^u(\mathfrak{t}) \cdot \mathbb{E}[\Lambda^u(\mathfrak{t},\theta)\Xi|\mathcal{F}_{\mathfrak{t}}]}{\Lambda^u(\mathfrak{t}) \cdot \mathbb{E}[\Lambda^u(\mathfrak{t},\theta)|\mathcal{F}_{\mathfrak{t}}]} = \mathbb{E}[\Lambda^u(\mathfrak{t},\theta)\Xi|\mathcal{F}_{\mathfrak{t}}]$$

(25)
$$= \mathbb{E}[\Lambda^v(\mathfrak{t},\theta)\Xi|\mathcal{F}_{\mathfrak{t}}] = \cdots = \mathbb{E}^v[\Xi|\mathcal{F}_{\mathfrak{t}}] \qquad \text{a.s.}$$

The second claim is proved similarly.

4.3. *Families directed downward.* For any given control strategy $v(\cdot) \in \mathfrak{U}$ and stopping rules $\mathfrak{t}, \theta$ with $0 \le \mathfrak{t} \le \theta \le T$, we shall denote by $\mathcal{V}_{[\mathfrak{t},\theta]}$ the set of admissible control strategies $u(\cdot)$ as in Lemma 4.3 (i.e., with $u(\cdot) = v(\cdot)$ a.e. on the stochastic interval $[[\mathfrak{t},\theta]]$).

We observe from (19), (8) and Lemma 4.3 that $Z^u(\theta)$ depends only on the values that the admissible control strategy $u(\cdot)$ takes over the stochastic interval $]]\theta, T]] := \{(s,\omega) \in [0,T] \times \Omega : \theta(\omega) < s \le T\}$ (its values over the stochastic interval $[[0,\theta]]$ are irrelevant for computing $Z^u(\theta)$). Thus, for *any* given admissible control strategy $v(\cdot) \in \mathfrak{U}$, we can write the upper value (11) of the game as

(26) $$\overline{V}(\theta) = \operatorname*{essinf}_{u \in \mathfrak{U}} Z^u(\theta) = \operatorname*{essinf}_{u \in \mathcal{V}_{[0,\theta]}} Z^u(\theta) \qquad \text{a.s.}$$

LEMMA 4.4.  *The family of random variables $\{Z^u(\theta)\}_{u \in \mathcal{V}_{[0,\theta]}}$ is* directed downward: *for any two $u^1(\cdot) \in \mathcal{V}_{[0,\theta]}$ and $u^2(\cdot) \in \mathcal{V}_{[0,\theta]}$, there exists an admissible control strategy $\widehat{u}(\cdot) \in \mathcal{V}_{[0,\theta]}$ such that we have a.s.*

$$Z^{\widehat{u}}(\theta) = Z^{u^1}(\theta) \wedge Z^{u^2}(\theta).$$

PROOF. Consider the event $A := \{Z^{u^1}(\theta) \le Z^{u^2}(\theta)\} \in \mathcal{F}_\theta$, and define an admissible control process $u(\cdot) \in \mathfrak{U}$ via $\widehat{u}(s,\omega) := v(s,\omega)$ for $0 \le s \le \theta(\omega)$,

(27)  $\widehat{u}(s,\omega) := u^1(s,\omega) \cdot 1_A(\omega) + u^2(s,\omega) \cdot 1_{A^c}(\omega)$ for $\theta(\omega) < s \le T$.

Consider also the stopping rule $\widehat{\tau}_\theta := \tau_\theta^{u^1} \cdot 1_A + \tau_\theta^{u^2} \cdot 1_{A^c} \in \mathcal{S}_{\theta,T}$ [notation of (22)]. Then from Lemma 4.3(ii) we have

$$Z^{\widehat{u}}(\theta) = \mathbb{E}^{\widehat{u}}[Y^{\widehat{u}}(\theta, \tau_\theta^{\widehat{u}})|\mathcal{F}_\theta]$$

$$= \mathbb{E}^{u^1}[Y^{u^1}(\theta, \tau_\theta^{\widehat{u}})|\mathcal{F}_\theta] \cdot 1_A + \mathbb{E}^{u^2}[Y^{u^2}(\theta, \tau_\theta^{\widehat{u}})|\mathcal{F}_\theta] \cdot 1_{A^c}$$

$$\le Z^{u^1}(\theta) \cdot 1_A + Z^{u^2}(\theta) \cdot 1_{A^c}$$

(28)
$$= \mathbb{E}^{u^1}[Y^{u^1}(\theta, \tau_\theta^{u^1})|\mathcal{F}_\theta] \cdot 1_A + \mathbb{E}^{u^2}[Y^{u^2}(\theta, \tau_\theta^{u^2})|\mathcal{F}_\theta] \cdot 1_{A^c}$$

$$= \mathbb{E}^{\widehat{u}}[Y^{\widehat{u}}(\theta, \tau_\theta^{u^1})|\mathcal{F}_\theta] \cdot 1_A + \mathbb{E}^{\widehat{u}}[Y^{\widehat{u}}(\theta, \tau_\theta^{u^2})|\mathcal{F}_\theta] \cdot 1_{A^c}$$

$$= \mathbb{E}^{\widehat{u}}[Y^{\widehat{u}}(\theta, \widehat{\tau}_\theta)|\mathcal{F}_\theta] \le Z^{\widehat{u}}(\theta),$$

thus, also $Z^{\widehat{u}}(\theta) = Z^{u^1}(\theta) \cdot 1_A + Z^{u^2}(\theta) \cdot 1_{A^c} = Z^{u^1}(\theta) \wedge Z^{u^2}(\theta)$, a.s. □

Now we can appeal to basic properties of the essential infimum [e.g., Neveu (1975), page 121], to obtain the following.

LEMMA 4.5. *For each $\theta \in \mathcal{S}$, there exists a sequence of admissible control processes $\{u^k(\cdot)\}_{k \in \mathbb{N}} \subset \mathcal{V}_{[0,\theta]}$, such that the corresponding sequence of random variables $\{Z^{u^k}(\theta)\}_{k \in \mathbb{N}}$ is decreasing, and the essential infimum in (26) becomes a limit:*

$$\overline{V}(\theta) = \lim_{k \to \infty} \downarrow Z^{u^k}(\theta) \qquad \text{a.s.} \tag{29}$$

**5. Existence and regularity of the game's value process.** For any given $\theta \in \mathcal{S}$, and with $\{u^k(\cdot)\}_{k \in \mathbb{N}} \subset \mathcal{V}_{[0,\theta]}$ the sequence of (29), let us look at the corresponding stopping rules

$$\tau_\theta^{u^k} := \inf\{s \in [\theta, T] : Z^{u^k}(s) = g(X(s))\}, \qquad k \in \mathbb{N},$$

via (22). Recall that we have $Z^{u^k}(\cdot) \geq Z^{u^\ell}(\cdot) \geq g(X(\cdot))$ for any integers $\ell \geq k$, thus, also $\tau_\theta^{u^k} \geq \tau_\theta^{u^\ell} \geq \theta$. In other words, the resulting sequence $\{\tau_\theta^{u^k}\}_{k \in \mathbb{N}}$ is decreasing, so the limit

$$\tau_\theta^* := \lim_{k \to \infty} \downarrow \tau_\theta^{u^k} \tag{30}$$

exists a.s. and defines a stopping rule in $\mathcal{S}_{\theta,T}$. The values of the process $u^k(\cdot)$ on the stochastic interval $[[0, \theta]]$ are irrelevant for computing $Z^{u^k}(s)$, $s \geq \theta$ or, for that matter, $\tau_\theta^{u^k}$. But clearly,

$$\tau_\theta^{u^k} = \inf\{s \in [\tau_\theta^*, T] : Z^{u^k}(s) = g(X(s))\}, \qquad k \in \mathbb{N},$$

holds a.s., so the values of $u^k(\cdot)$ on $[[0, \tau_\theta^*]]$ are irrelevant for computing $\tau_\theta^{u^k}$, $k \in \mathbb{N}$.

Thus, there exists a sequence $\{u^k(\cdot)\}_{k \in \mathbb{N}} \subset \mathcal{V}_{[0,\tau_\theta^*]}$ of admissible control strategies, which agree with the given control strategy $v(\cdot) \in \mathfrak{U}$ on the stochastic interval $[[0, \tau_\theta^*]]$, and for which (30) holds.

We are ready to state and prove our first result.

THEOREM 5.1. *The game has a value: for every $\theta \in \mathcal{S}$, we have $\overline{V}(\theta) = \underline{V}(\theta)$, a.s. In particular, $\overline{V} = \underline{V}$ in (9). A bit more generally, for every $\mathfrak{t} \in \mathcal{S}$ and any $\theta \in \mathcal{S}_{\mathfrak{t},T}$, we have, almost surely,*

$$\operatorname*{essinf}_{u \in \mathfrak{U}} \operatorname*{esssup}_{\tau \in \mathcal{S}_{\theta,T}} \mathbb{E}^u(Y^u(\mathfrak{t}, \tau) | \mathcal{F}_\mathfrak{t}) = \operatorname*{esssup}_{\tau \in \mathcal{S}_{\theta,T}} \operatorname*{essinf}_{u \in \mathfrak{U}} \mathbb{E}^u(Y^u(\mathfrak{t}, \tau) | \mathcal{F}_\mathfrak{t}). \tag{31}$$

PROOF. From the preceding remarks, we get the a.s. comparisons

$$\overline{V}(\theta) \leq \mathbb{E}^{u^k}[Y^{u^k}(\theta, \tau_\theta^{u^k})|\mathcal{F}_\theta]$$
$$= \mathbb{E}[\Lambda^{u^k}(\theta, \tau_\theta^{u^k})Y^{u^k}(\theta, \tau_\theta^{u^k})|\mathcal{F}_\theta]$$
$$= \mathbb{E}\left[\Lambda^v(\theta, \tau_\theta^*)\Lambda^{u^k}(\tau_\theta^*, \tau_\theta^{u^k})\left\{Y^v(\theta, \tau_\theta^*) + \int_{\tau_\theta^*}^{\tau_\theta^{u^k}} h(s, X, u_s^k)\,ds\right\}\bigg|\mathcal{F}_\theta\right]$$

for every $k \in \mathbb{N}$; recall the computation (25). Passing to the limit as $k \to \infty$, and using (30), the boundedness of $\sigma^{-1}$, $f$, $h$, and the dominated convergence theorem, we obtain the a.s. comparisons

$$\overline{V}(\theta) \leq \mathbb{E}[\Lambda^v(\theta, \tau_\theta^*)Y^v(\theta, \tau_\theta^*)|\mathcal{F}_\theta] = \mathbb{E}^v[Y^v(\theta, \tau_\theta^*)|\mathcal{F}_\theta].$$

Because $v(\cdot)$ is arbitrary, we can take the infimum of the right-hand side of this inequality over $v(\cdot) \in \mathfrak{U}$, and conclude

$$\overline{V}(\theta) \leq \operatorname*{essinf}_{v \in \mathfrak{U}} \mathbb{E}^v[Y^v(\theta, \tau_\theta^*)|\mathcal{F}_\theta]$$
$$\leq \operatorname*{esssup}_{\tau \in \mathcal{S}_{\theta,T}} \operatorname*{essinf}_{v \in \mathfrak{U}} \mathbb{E}^v[Y^v(\theta, \tau)|\mathcal{F}_\theta] = \underline{V}(\theta).$$

The reverse inequality $\overline{V}(\theta) \geq \underline{V}(\theta)$ is obvious, so we obtain the first claim of the theorem, namely, $\overline{V}(\theta) = \underline{V}(\theta)$ a.s.

• As for (31), let us observe that for every given $u(\cdot) \in \mathfrak{U}$ we have, on the strength of Proposition 4.2, the a.s. comparisons

$$\operatorname*{essinf}_{w \in \mathfrak{U}} \operatorname*{esssup}_{\tau \in \mathcal{S}_{\theta,T}} \mathbb{E}^w\left(Y^w(\theta, \tau) + \int_{\mathfrak{t}}^\theta h(s, X, w_s)\,ds\bigg|\mathcal{F}_{\mathfrak{t}}\right)$$
$$\leq \operatorname*{esssup}_{\tau \in \mathcal{S}_{\theta,T}} \mathbb{E}^u\left(Y^u(\theta, \tau) + \int_{\mathfrak{t}}^\theta h(s, X, u_s)\,ds\bigg|\mathcal{F}_{\mathfrak{t}}\right)$$
$$\leq \mathbb{E}^u\left(\operatorname*{esssup}_{\tau \in \mathcal{S}_{\theta,T}} \mathbb{E}^u[Y^u(\theta, \tau)|\mathcal{F}_\theta] + \int_{\mathfrak{t}}^\theta h(s, X, u_s)\,ds\bigg|\mathcal{F}_{\mathfrak{t}}\right)$$
$$= \mathbb{E}^u\left(\mathbb{E}^u[Y^u(\theta, \tau_\theta^u)|\mathcal{F}_\theta] + \int_{\mathfrak{t}}^\theta h(s, X, u_s)\,ds\bigg|\mathcal{F}_{\mathfrak{t}}\right)$$
$$= \mathbb{E}^u\left(Y^u(\theta, \tau_\theta^u) + \int_{\mathfrak{t}}^\theta h(s, X, u_s)\,ds\bigg|\mathcal{F}_{\mathfrak{t}}\right)$$
$$= \mathbb{E}^u[Y^u(\mathfrak{t}, \tau_\theta^u)|\mathcal{F}_{\mathfrak{t}}].$$

Now repeat the previous argument: fix $v(\cdot) \in \mathfrak{U}$, write this inequality with $u(\cdot)$ replaced by $u^k(\cdot) \in \mathcal{V}_{[0,\tau_\theta^*]}$ [the sequence of (29), (30)] for every $k \in \mathbb{N}$, and observe that the last term in the above string is now equal to $\mathbb{E}^v[Y^v(\mathfrak{t}, \tau_\theta^{u^k})|\mathcal{F}_\mathfrak{t}]$. Then pass to the limit as $k \to \infty$ to get, a.s.,

$$\operatorname*{essinf}_{w \in \mathfrak{U}} \operatorname*{esssup}_{\tau \in \mathcal{S}_{\theta,T}} \mathbb{E}^w\left(Y^w(\theta, \tau) + \int_\mathfrak{t}^\theta h(s, X, w_s)\,ds \Big| \mathcal{F}_\mathfrak{t}\right) \leq \mathbb{E}^v[Y^v(\mathfrak{t}, \tau_\theta^*)|\mathcal{F}_\mathfrak{t}].$$

The arbitrariness of $v(\cdot)$ allows us to take the (essential) infimum of the right-hand side over $v(\cdot) \in \mathfrak{U}$, and obtain

$$\operatorname*{essinf}_{w \in \mathfrak{U}} \operatorname*{esssup}_{\tau \in \mathcal{S}_{\theta,T}} \mathbb{E}^w[Y^w(\mathfrak{t}, \tau)|\mathcal{F}_\mathfrak{t}] \leq \operatorname*{essinf}_{v \in \mathfrak{U}} \mathbb{E}^v[Y^v(\mathfrak{t}, \tau_\theta^*)|\mathcal{F}_\mathfrak{t}]$$

$$\leq \operatorname*{esssup}_{\tau \in \mathcal{S}_{\theta,T}} \operatorname*{essinf}_{v \in \mathfrak{U}} \mathbb{E}^v[Y^v(\mathfrak{t}, \tau)|\mathcal{F}_\mathfrak{t}],$$

that is, the inequality ($\leq$) of (31); once again, the reverse inequality is obvious. □

From now on we shall denote by $V(\cdot) = \overline{V}(\cdot) = \underline{V}(\cdot)$ the common value process of this game, and write $V = V(0)$.

PROPOSITION 5.2. *The value process $V(\cdot)$ is right-continuous.*

PROOF. The Snell Envelope $Q^u(\cdot)$ of (20) can be taken in its RCLL modification, as we have already done; so the same is the case for the process $Z^u(\cdot)$ of (19). Consequently, we obtain $\limsup_{s \downarrow t} V(s) \leq \lim_{s \downarrow t} Z^u(s) = Z^u(t)$, a.s. Taking the infimum over $u(\cdot) \in \mathfrak{U}$, we deduce $\limsup_{s \downarrow t} V(s) \leq V(t)$, a.s.

In order to show that the reverse inequality

(32) $$\liminf_{s \downarrow t} V(s) \geq V(t) \qquad \text{a.s.,}$$

also holds, recall the submartingale property of (17) and (18) and deduce from it, and from Proposition 1.3.14 in Karatzas and Shreve (1991), that the right-hand limits

$$J(t+, \tau) := \lim_{s \downarrow t} J(s, \tau) \quad \text{on } \{t < \tau\}, \qquad J(t+, \tau) := g(X(\tau)) \quad \text{on } \{t = \tau\}$$

exist and are finite, a.s. on the respective events. Now for any $t \in [0, T]$ and every stopping rule $\tau \in \mathcal{S}_{t,T}$, recall (16) to obtain

$$\liminf_{s \downarrow t} V(s) \geq \liminf_{s \downarrow t} J(s, s \vee \tau)$$

$$= \liminf_{s \downarrow t} J(s, \tau) \cdot 1_{\{t < \tau\}} + \liminf_{s \downarrow t} J(s, s) \cdot 1_{\{t = \tau\}}.$$

But on the event $\{t = \tau\}$, we have almost surely

$$\liminf_{s \downarrow t} J(s, s) = \liminf_{s \downarrow t} g(X(s)) = \lim_{s \downarrow t} g(X(s)) = g(X(t)) = J(t, t)$$

by the continuity of $g(\cdot)$; whereas on the event $\{t < \tau\}$, we have the a.s. equalities $\liminf_{s \downarrow t} J(s, \tau) = \lim_{s \downarrow t} J(s, \tau) = J(t+, \tau)$. Recalling (18), we obtain from the bounded convergence theorem the a.s. comparisons

$$\liminf_{s \downarrow t} V(s) \geq \lim_{s \downarrow t} J(s, \tau) = \mathbb{E}^u\left(\lim_{s \downarrow t} J(s, \tau) \Big| \mathcal{F}_{t+}\right) = \mathbb{E}^u\left(\lim_{s \downarrow t} J(s, \tau) \Big| \mathcal{F}_t\right)$$

$$= \mathbb{E}^u\left[\lim_{s \downarrow t}\left(J(s, \tau) + \int_t^s h(r, X, u_r)\,dr\right) \Big| \mathcal{F}_t\right]$$

$$= \lim_{s \downarrow t} \mathbb{E}^u\left[J(s, \tau) + \int_t^s h(r, X, u_r)\,dr \Big| \mathcal{F}_t\right] \geq J(t, \tau).$$

We have used here the right-continuity of the augmented Brownian filtration [Karatzas and Shreve (1991), pages 89–92]. The stopping rule $\tau \in \mathcal{S}_{t,T}$ is arbitrary in these comparisons; taking the (essential) supremum over $\mathcal{S}_{t,T}$ and recalling (16), we arrive at the desired inequality (32). □

5.1. *Some elementary submartingales.* By analogy with (22), let us introduce now for each $\mathfrak{t} \in \mathcal{S}$ and $0 \leq \varepsilon < 1$ the stopping rules

(33) $\quad \varrho_\mathfrak{t}(\varepsilon) := \inf\{s \in [\mathfrak{t}, T] : g(X(s)) \geq V(s) - \varepsilon\}, \qquad \varrho_\mathfrak{t} := \varrho_\mathfrak{t}(0).$

Since

(34) $$V(\cdot) = \operatorname*{essinf}_{u \in \mathfrak{U}} Z^u(\cdot) \geq g(X(\cdot))$$

holds a.s. thanks to (26), we have also

(35) $\quad \varrho_\mathfrak{t} \vee \tau_\mathfrak{t}^u(\varepsilon) \leq \tau_\mathfrak{t}^u, \qquad \varrho_\mathfrak{t}(\varepsilon) \leq \tau_\mathfrak{t}^u(\varepsilon) \wedge \varrho_\mathfrak{t}.$

For each admissible control strategy $u(\cdot) \in \mathfrak{U}$, let us introduce the family of random variables

(36) $\quad R^u(\mathfrak{t}) := V(\mathfrak{t}) + \int_0^\mathfrak{t} h(s, X, u_s)\,ds \geq Y^u(\mathfrak{t}), \qquad \mathfrak{t} \in \mathcal{S}.$

For any time $\mathfrak{t} \in \mathcal{S}$, the quantity $R^u(\mathfrak{t})$ represents the cumulative cost to the controller of using the strategy $u(\cdot)$ on $[[0, \mathfrak{t}]]$, plus the game's value at that time.

PROPOSITION 5.3. *For each $u(\cdot) \in \mathfrak{U}$, the process $R^u(\cdot \wedge \varrho_0)$ is a $\mathbb{P}^u$-submartingale. A bit more generally, for any stopping rules $\mathfrak{t}, \vartheta$ with $\mathfrak{t} \leq \vartheta \leq \varrho_\mathfrak{t}$, we have*

(37) $\quad \mathbb{E}^u[R^u(\vartheta)|\mathcal{F}_\mathfrak{t}] \geq R^u(\mathfrak{t}) \qquad a.s.,$

*or, equivalently,*

(38) $\quad \mathbb{E}^u\left[V(\vartheta) + \int_\mathfrak{t}^\vartheta h(s, X, u_s)\,ds \Big| \mathcal{F}_\mathfrak{t}\right] \geq V(\mathfrak{t}) \qquad a.s.$

*Furthermore, for any stopping rules* $\mathfrak{s}, \mathfrak{t}, \vartheta$ *with* $0 \le \mathfrak{s} \le \mathfrak{t} \le \vartheta \le \varrho_\mathfrak{t}$, *we have almost surely*

$$(39) \quad \operatorname*{essinf}_{u \in \mathfrak{U}} \mathbb{E}^u \bigg[ V(\vartheta) + \int_\mathfrak{s}^\vartheta h(s, X, u_s) \, ds \bigg| \mathcal{F}_\mathfrak{s} \bigg]$$
$$\ge \operatorname*{essinf}_{u \in \mathfrak{U}} \mathbb{E}^u \bigg[ V(\mathfrak{t}) + \int_\mathfrak{s}^\mathfrak{t} h(s, X, u_s) \, ds \bigg| \mathcal{F}_\mathfrak{s} \bigg].$$

PROOF. For any admissible control strategy $u(\cdot) \in \mathfrak{U}$, and for any stopping rules $\mathfrak{t}, \vartheta$ with $0 \le \mathfrak{t} \le \vartheta \le \varrho_\mathfrak{t}$, we have $\mathbb{E}^u[Q^u(\vartheta)|\mathcal{F}_\mathfrak{t}] = Q^u(\mathfrak{t})$ or, equivalently,

$$(40) \quad \mathbb{E}^u \bigg[ Z^u(\vartheta) + \int_\mathfrak{t}^\vartheta h(s, X, u_s) \, ds \bigg| \mathcal{F}_\mathfrak{t} \bigg] = Z^u(\mathfrak{t}) \ge V(\mathfrak{t}) \qquad \text{a.s.}$$

from (21), (35) and Propositions 4.1 and 4.2. Now *fix* a control strategy $v(\cdot) \in \mathfrak{U}$, and denote again by $\mathcal{V}_{[\mathfrak{t},\vartheta]}$ the set of admissible control strategies $u(\cdot)$ as in Lemma 4.3 that agree with it [i.e., satisfy $u(\cdot) = v(\cdot)$ a.e.] on the stochastic interval $[[\mathfrak{t}, \vartheta]]$. From this result and (40), we obtain

$$(41) \quad \mathbb{E}^v \bigg[ Z^u(\vartheta) + \int_\mathfrak{t}^\vartheta h(s, X, v_s) \, ds \bigg| \mathcal{F}_\mathfrak{t} \bigg] = Z^u(\mathfrak{t}) \ge V(\mathfrak{t}) \qquad \text{a.s.}$$

Now select some sequence $\{u^k(\cdot)\}_{k \in \mathbb{N}} \subset \mathcal{V}_{[\mathfrak{t},\vartheta]}$ as in (29) of Lemma 4.5, substitute $u^k(\cdot)$ for $u(\cdot)$ in (41), let $k \to \infty$, and appeal to the bounded convergence theorem for conditional expectations to obtain

$$\mathbb{E}^v \bigg[ V(\vartheta) + \int_\mathfrak{t}^\vartheta h(s, X, v_s) \, ds \bigg| \mathcal{F}_\mathfrak{t} \bigg] \ge V(\mathfrak{t}) \qquad \text{a.s.}$$

This gives (38), therefore, also

$$\mathbb{E}^u \bigg[ V(\vartheta) + \int_\mathfrak{s}^\vartheta h(s, X, u_s) \, ds \bigg| \mathcal{F}_\mathfrak{s} \bigg] \ge \mathbb{E}^u \bigg[ V(\mathfrak{t}) + \int_\mathfrak{s}^\mathfrak{t} h(s, X, u_s) \, ds \bigg| \mathcal{F}_\mathfrak{s} \bigg],$$

for all $u(\cdot) \in \mathfrak{U}$. The claim (39) follows now by taking essential infima over $u(\cdot) \in \mathfrak{U}$ on both sides. $\square$

PROPOSITION 5.4. *For every* $\mathfrak{t} \in \mathcal{S}$, *we have*

$$(42) \quad V(\mathfrak{t}) = \operatorname*{essinf}_{u \in \mathfrak{U}} \mathbb{E}^u \bigg( g(X(\varrho_\mathfrak{t})) + \int_\mathfrak{t}^{\varrho_\mathfrak{t}} h(s, X, u_s) \, ds \bigg| \mathcal{F}_\mathfrak{t} \bigg) \qquad \text{a.s.}$$

*As a consequence,*

$$(43) \quad V(\mathfrak{t}) = \operatorname*{essinf}_{u \in \mathfrak{U}} \mathbb{E}^u \bigg( V(\varrho_\mathfrak{t}) + \int_\mathfrak{t}^{\varrho_\mathfrak{t}} h(s, X, u_s) \, ds \bigg| \mathcal{F}_\mathfrak{t} \bigg) \qquad \text{a.s.}$$

*and for any given* $v(\cdot) \in \mathfrak{U}$, *we get in the notation of* (26):

$$(44) \quad R^v(\mathfrak{t}) = \operatorname*{essinf}_{u \in \mathcal{V}_{[0,\mathfrak{t}]}} \mathbb{E}^u(R^u(\varrho_\mathfrak{t})|\mathcal{F}_\mathfrak{t}) \qquad \text{a.s.}$$

PROOF. The definition (11) for the upper value of the game gives the inequality ($\geq$) in (42). For the reverse inequality ($\leq$), write (38) of Proposition 5.3 with $\vartheta = \varrho_{\mathfrak{t}}$ and recall the a.s. equality $V(\varrho_{\mathfrak{t}}) = g(X(\varrho_{\mathfrak{t}}))$, a consequence of the definition of $\varrho_{\mathfrak{t}}$ in (33) and the right-continuity of $V(\cdot)$ from Proposition 5.2; the result is

$$V(\mathfrak{t}) \leq \mathbb{E}^u \bigg( V(\varrho_{\mathfrak{t}}) + \int_{\mathfrak{t}}^{\varrho_{\mathfrak{t}}} h(s, X, u_s)\, ds \bigg| \mathcal{F}_{\mathfrak{t}} \bigg)$$
$$= \mathbb{E}^u \bigg( g(X(\varrho_{\mathfrak{t}})) + \int_{\mathfrak{t}}^{\varrho_{\mathfrak{t}}} h(s, X, u_s)\, ds \bigg| \mathcal{F}_{\mathfrak{t}} \bigg) \qquad \text{a.s.}$$

for every $u(\cdot) \in \mathfrak{U}$. Now (42) and (43) follow directly, and so does (44). $\square$

REMARK 5.1. Proposition 5.3 implies that the process $R^u(\cdot \wedge \varrho_0)$, which is right-continuous by virtue of Proposition 5.2, admits left-limits on $(0, T]$ almost surely; cf. Proposition 1.3.14 in Karatzas and Shreve (1991). Thus, the process $R^u(\cdot \wedge \varrho_0)$ is a $\mathbb{P}^u$-submartingale with RCLL paths, and the process $V(\cdot \wedge \varrho_0)$ has RCLL paths as well.

**6. Some properties of the value process.** We shall derive in this section some further properties of $V(\cdot)$, the value process of the stochastic game. These will be crucial in characterizing, then constructing, a saddle point for the game in Sections 7 and 9, respectively.

Our first such result provides inequalities in the reverse direction of those in (37) and (38), but for more general stopping rules and with appropriate modifications.

PROPOSITION 6.1. *For any stopping rules $\mathfrak{t}, \theta$ with $0 \leq \mathfrak{t} \leq \theta \leq T$, and any admissible control process $u(\cdot) \in \mathfrak{U}$, we have*

$$(45) \qquad \mathbb{E}^u[R^u(\theta)|\mathcal{F}_{\mathfrak{t}}] \leq \operatorname*{esssup}_{\tau \in \mathcal{S}_{\mathfrak{t}, T}} \mathbb{E}^u(Y^u(\tau)|\mathcal{F}_{\mathfrak{t}})$$

*and*

$$(46) \quad \mathbb{E}^u \bigg[ V(\theta) + \int_{\mathfrak{t}}^{\theta} h(s, X, u_s)\, ds \bigg| \mathcal{F}_{\mathfrak{t}} \bigg] \leq \operatorname*{esssup}_{\tau \in \mathcal{S}_{\mathfrak{t}, T}} \mathbb{E}^u(Y^u(\mathfrak{t}, \tau)|\mathcal{F}_{\mathfrak{t}}) = Z^u(\mathfrak{t})$$

*almost surely. We also have*

$$(47) \qquad \operatorname*{essinf}_{u \in \mathfrak{U}} \mathbb{E}^u \bigg[ V(\theta) + \int_{\mathfrak{t}}^{\theta} h(s, X, u_s)\, ds \bigg| \mathcal{F}_{\mathfrak{t}} \bigg] \leq V(\mathfrak{t}) \qquad a.s.$$

*and*

$$(48) \qquad \operatorname*{essinf}_{u \in \mathcal{V}_{[0,\mathfrak{t}]}} \mathbb{E}^u \bigg[ V(\theta) + \int_{\mathfrak{t}}^{\theta} h(s, X, u_s)\, ds \bigg| \mathcal{F}_{\mathfrak{t}} \bigg] \leq V(\mathfrak{t}) \qquad a.s.$$

*for any given $v(\cdot) \in \mathfrak{U}$ in the notation used in (26).*

PROOF. We recall from (26), (19) and Theorem 5.1 that $V(\theta) = \operatorname{ess\,inf}_{u \in \mathfrak{U}} Z^u(\theta)$; and from Proposition 4.1 that, for any given $u(\cdot) \in \mathfrak{U}$, the process $Q^u(\cdot) = Z^u(\cdot) + \int_0^{\cdot} h(s, X, u_s) \, ds$ is a $\mathbb{P}^u$-supermartingale. We have, therefore,

$$\begin{aligned}
\mathbb{E}^u[R^u(\theta)|\mathcal{F}_{\mathfrak{t}}] &= \mathbb{E}^u\left[V(\theta) + \int_0^{\theta} h(s, X, u_s) \, ds \Big| \mathcal{F}_{\mathfrak{t}}\right] \\
&\leq \mathbb{E}^u\left[Z^u(\theta) + \int_0^{\theta} h(s, X, u_s) \, ds \Big| \mathcal{F}_{\mathfrak{t}}\right] \\
&\leq Z^u(\mathfrak{t}) + \int_0^{\mathfrak{t}} h(s, X, u_s) \, ds \\
&= \operatorname*{ess\,sup}_{\tau \in \mathcal{S}_{\mathfrak{t},T}} \mathbb{E}^u(Y^u(\tau)|\mathcal{F}_{\mathfrak{t}}),
\end{aligned}$$
(49)

which is (45). Now (46) is a direct consequence; and (47) and (48) follow by taking essential infima over $u(\cdot)$ in $\mathfrak{U}$ and in $\mathcal{V}_{[0,\mathfrak{t}]}$, respectively. □

We have also the following result, which supplements the "value identity" of equation (31). In this equation the common value is at most $V(\mathfrak{t})$, as we are taking supremum over a class of stopping rules, $\mathcal{S}_{\theta,T}$, which is smaller than the class $\mathcal{S}_{\mathfrak{t},T}$ appearing in (11) and (12). The next result tells us exactly how smaller than $V(\mathfrak{t})$ this common value is: it is given by the left-hand side of (47).

PROPOSITION 6.2. *For any stopping rules* $\mathfrak{t}, \theta$ *with* $0 \leq \mathfrak{t} \leq \theta \leq T$, *we have almost surely*

$$\begin{aligned}
&\operatorname*{ess\,inf}_{u \in \mathfrak{U}} \mathbb{E}^u\left[V(\theta) + \int_{\mathfrak{t}}^{\theta} h(s, X, u_s) \, ds \Big| \mathcal{F}_{\mathfrak{t}}\right] \\
&= \operatorname*{ess\,inf}_{u \in \mathfrak{U}} \operatorname*{ess\,sup}_{\tau \in \mathcal{S}_{\theta,T}} \mathbb{E}^u(Y^u(\mathfrak{t},\tau)|\mathcal{F}_{\mathfrak{t}}) \\
&= \operatorname*{ess\,sup}_{\tau \in \mathcal{S}_{\theta,T}} \operatorname*{ess\,inf}_{u \in \mathfrak{U}} \mathbb{E}^u(Y^u(\mathfrak{t},\tau)|\mathcal{F}_{\mathfrak{t}}).
\end{aligned}$$
(50)

PROOF. The second equality is, of course, that of (31). For the first, note that Proposition 5.4 gives $V(\theta) \leq \mathbb{E}^u(g(X(\varrho_\theta)) + \int_\theta^{\varrho_\theta} h(s, X, u_s) \, ds | \mathcal{F}_\theta)$ a.s., for every admissible control strategy $u(\cdot) \in \mathfrak{U}$, thus, also

$$\begin{aligned}
&\mathbb{E}^u\left[V(\theta) + \int_{\mathfrak{t}}^{\theta} h(s, X, u_s) \, ds \Big| \mathcal{F}_{\mathfrak{t}}\right] \\
&\leq \mathbb{E}^u\left(g(X(\varrho_\theta)) + \int_{\mathfrak{t}}^{\varrho_\theta} h(s, X, u_s) \, ds \Big| \mathcal{F}_{\mathfrak{t}}\right)
\end{aligned}$$
(51)

$$= \mathbb{E}^u(Y^u(\mathfrak{t}, \varrho_\theta)|\mathcal{F}_\mathfrak{t})$$
$$\leq \operatorname*{esssup}_{\tau \in \mathcal{S}_{\theta,T}} \mathbb{E}^u(Y^u(\mathfrak{t}, \tau)|\mathcal{F}_\mathfrak{t}) \qquad \text{a.s.}$$

Taking essential infima on both sides over $u(\cdot) \in \mathfrak{U}$, we arrive at the inequality ($\leq$) in (50).

For the reverse inequality, note from the definition of (19) that

$$Z^u(\theta) + \int_\mathfrak{t}^\theta h(s, X, u_s)\,ds \geq \mathbb{E}^u(Y^u(\mathfrak{t}, \tau)|\mathcal{F}_\theta)$$

holds a.s. for every $u(\cdot) \in \mathfrak{U}$ and every $\tau \in \mathcal{S}_{\theta,T}$ [in fact, with equality for the stopping rule, $\tau = \tau_\theta^u$ of (22)]. Taking conditional expectations with respect to $\mathcal{F}_\mathfrak{t}$ on both sides, we obtain

$$(52) \quad \mathbb{E}^u\left(Z^u(\theta) + \int_\mathfrak{t}^\theta h(s, X, u_s)\,ds \Big|\mathcal{F}_\mathfrak{t}\right) \geq \mathbb{E}^u(Y^u(\mathfrak{t}, \tau)|\mathcal{F}_\mathfrak{t}) \qquad \text{a.s.}$$

for all $\tau \in \mathcal{S}_{\theta,T}$, again with equality for $\tau = \tau_\theta^u$; therefore,

$$(53) \quad \mathbb{E}^u\left(Z^u(\theta) + \int_\mathfrak{t}^\theta h(s, X, u_s)\,ds \Big|\mathcal{F}_\mathfrak{t}\right) = \operatorname*{esssup}_{\tau \in \mathcal{S}_{\theta,T}} \mathbb{E}^u(Y^u(\mathfrak{t}, \tau)|\mathcal{F}_\mathfrak{t}) \qquad \text{a.s.}$$

Fix now an admissible control strategy $v(\cdot) \in \mathfrak{U}$, and consider a sequence $\{u^k(\cdot)\}_{k \in \mathbb{N}} \subset \mathcal{V}_{[\mathfrak{t},\theta]}$ such that $V(\theta) = \lim_{k \to \infty} \downarrow Z^{u^k}(\theta)$ a.s., in the manner of (29) in Lemma 4.5. Write (53) with $u^k(\cdot)$ in place of $u(\cdot)$ and recall property (24) of Lemma 4.3 to obtain

$$\mathbb{E}^v\left(Z^{u^k}(\theta) + \int_\mathfrak{t}^\theta h(s, X, v_s)\,ds \Big|\mathcal{F}_\mathfrak{t}\right)$$
$$= \mathbb{E}^{u^k}\left(Z^{u^k}(\theta) + \int_\mathfrak{t}^\theta h(s, X, u_s^k)\,ds \Big|\mathcal{F}_\mathfrak{t}\right)$$
$$= \operatorname*{esssup}_{\tau \in \mathcal{S}_{\theta,T}} \mathbb{E}^{u^k}(Y^{u^k}(\mathfrak{t}, \tau)|\mathcal{F}_\mathfrak{t})$$
$$\geq \operatorname*{essinf}_{u \in \mathfrak{U}} \operatorname*{esssup}_{\tau \in \mathcal{S}_{\theta,T}} \mathbb{E}^u(Y^u(\mathfrak{t}, \tau)|\mathcal{F}_\mathfrak{t}) \qquad \text{a.s.}$$

for every $k \in \mathbb{N}$. Now let $k \to \infty$ and use the bounded convergence theorem, to obtain

$$\mathbb{E}^v\left(V(\theta) + \int_\mathfrak{t}^\theta h(s, X, v_s)\,ds \Big|\mathcal{F}_\mathfrak{t}\right) \geq \operatorname*{essinf}_{u \in \mathfrak{U}} \operatorname*{esssup}_{\tau \in \mathcal{S}_{\theta,T}} \mathbb{E}^u(Y^u(\mathfrak{t}, \tau)|\mathcal{F}_\mathfrak{t}) \qquad \text{a.s.}$$

Since $v(\cdot) \in \mathfrak{U}$ is an arbitrary control strategy, all that remains at this point is to take the essential infimum of the left-hand side with respect to $v(\cdot) \in \mathfrak{U}$, and we are done. □

We are ready for the main result of this section. It says that $\inf_{u\in\mathfrak{U}} \mathbb{E}^u(R^u(\cdot))$, the best that the controller can achieve in terms of minimizing expected "running cost plus current value," does not increase with time; at best, this quantity is "flat up to $\varrho_0$," the first time the game's value equals the reward obtained by terminating the game.

THEOREM 6.3. *For any stopping rules $\mathfrak{t}, \theta$ with $0 \leq \mathfrak{t} \leq \theta \leq T$, we have*

$$\operatorname*{essinf}_{u\in\mathfrak{U}} \mathbb{E}^u(R^u(\theta)|\mathcal{F}_\mathfrak{t}) \leq R^v(\mathfrak{t}) \qquad \text{a.s.} \tag{54}$$

*for any $v(\cdot) \in \mathfrak{U}$, as well as*

$$\inf_{u\in\mathfrak{U}} \mathbb{E}^u(R^u(\theta)) \leq \inf_{u\in\mathfrak{U}} \mathbb{E}^u(R^u(\mathfrak{t})) \leq V(0). \tag{55}$$

*The first (resp., the second) of the inequalities in (55) is valid as equality if $\theta \leq \varrho_\mathfrak{t}$ (resp., $\mathfrak{t} \leq \varrho_0$) also holds.*

*A bit more generally, for any stopping rules $\mathfrak{s}, \mathfrak{t}, \theta$ with $0 \leq \mathfrak{s} \leq \mathfrak{t} \leq \theta \leq T$, we have the a.s. comparisons*

$$\begin{aligned}
\operatorname*{essinf}_{u\in\mathfrak{U}} \mathbb{E}^u &\left[ V(\theta) + \int_\mathfrak{s}^\theta h(s, X, u_s)\, ds \Big| \mathcal{F}_\mathfrak{s} \right] \\
&\leq \operatorname*{essinf}_{u\in\mathfrak{U}} \mathbb{E}^u \left[ V(\mathfrak{t}) + \int_\mathfrak{s}^\mathfrak{t} h(s, X, u_s)\, ds \Big| \mathcal{F}_\mathfrak{s} \right] \leq V(\mathfrak{s}).
\end{aligned} \tag{56}$$

*The first (resp., the second) of the inequalities in (56) is valid as an equality on the event $\{\theta \leq \varrho_\mathfrak{t}\}$ (resp., $\{\mathfrak{t} \leq \varrho_\mathfrak{s}\}$).*

PROOF. With $v(\cdot) \in \mathfrak{U}$ fixed, and with $\mathcal{V}_{[0,\mathfrak{t}]}$ as in Lemma 4.5, we have

$$\begin{aligned}
\operatorname*{essinf}_{u\in\mathfrak{U}} &\mathbb{E}^u(R^u(\theta)|\mathcal{F}_\mathfrak{t}) \\
&= \operatorname*{essinf}_{u\in\mathfrak{U}} \left( \mathbb{E}^u\left[V(\theta) + \int_\mathfrak{t}^\theta h(s, X, u_s)\,ds \Big| \mathcal{F}_\mathfrak{t}\right] + \int_0^\mathfrak{t} h(s, X, u_s)\,ds \right) \\
&\leq \operatorname*{essinf}_{u\in\mathcal{V}_{[0,\mathfrak{t}]}} \mathbb{E}^u\left[V(\theta) + \int_\mathfrak{t}^\theta h(s, X, u_s)\,ds \Big| \mathcal{F}_\mathfrak{t}\right] + \int_0^\mathfrak{t} h(s, X, v_s)\,ds \\
&\leq V(\mathfrak{t}) + \int_0^\mathfrak{t} h(s, X, v_s)\,ds = R^v(\mathfrak{t}) \qquad \text{a.s.}
\end{aligned}$$

where the penultimate comparison comes from (48). This proves (54).

To obtain the first inequality in (56), observe that (49) gives

$$\mathbb{E}^u\left[V(\theta) + \int_\mathfrak{s}^\theta h(s, X, u_s)\,ds \Big| \mathcal{F}_\mathfrak{s}\right] \leq \mathbb{E}^u\left[Z^u(\mathfrak{t}) + \int_\mathfrak{s}^\mathfrak{t} h(s, X, u_s)\,ds \Big| \mathcal{F}_\mathfrak{s}\right] \qquad \text{a.s.}$$

for all $u(\cdot) \in \mathfrak{U}$. Proceeding just as before, with $v(\cdot) \in \mathfrak{U}$ arbitrary but fixed, and with a sequence $\{u^k(\cdot)\}_{k \in \mathbb{N}} \subset \mathcal{V}_{[0,\mathfrak{t}]}$ such that $V(\mathfrak{t}) = \lim_{k \to \infty} \downarrow Z^{u^k}(\mathfrak{t})$ holds almost surely, as in Lemma 4.5, we have

$$\operatorname*{essinf}_{u \in \mathfrak{U}} \mathbb{E}^u \left[ V(\theta) + \int_{\mathfrak{s}}^{\theta} h(s, X, u_s) \, ds \Big| \mathcal{F}_{\mathfrak{s}} \right]$$
$$\leq \mathbb{E}^{u^k} \left[ V(\theta) + \int_{\mathfrak{s}}^{\theta} h(s, X, u_s^k) \, ds \Big| \mathcal{F}_{\mathfrak{s}} \right]$$
$$\leq \mathbb{E}^{u^k} \left[ Z^{u^k}(\mathfrak{t}) + \int_{\mathfrak{s}}^{\mathfrak{t}} h(s, X, u_s^k) \, ds \Big| \mathcal{F}_{\mathfrak{s}} \right]$$
$$= \mathbb{E}^v \left[ Z^{u^k}(\mathfrak{t}) + \int_{\mathfrak{s}}^{\mathfrak{t}} h(s, X, v_s) \, ds \Big| \mathcal{F}_{\mathfrak{s}} \right]$$

for every $k \in \mathbb{N}$, thus, also

$$\operatorname*{essinf}_{u \in \mathfrak{U}} \mathbb{E}^u \left[ V(\theta) + \int_{\mathfrak{s}}^{\theta} h(s, X, u_s) \, ds \Big| \mathcal{F}_{\mathfrak{s}} \right] \leq \mathbb{E}^v \left[ V(\mathfrak{t}) + \int_{\mathfrak{s}}^{\mathfrak{t}} h(s, X, v_s) \, ds \Big| \mathcal{F}_{\mathfrak{s}} \right]$$

in the limit as $k \to \infty$. Take the essential infimum of the right-hand side over $v(\cdot) \in \mathfrak{U}$ to obtain the desired a.s. inequality

$$\operatorname*{essinf}_{u \in \mathfrak{U}} \mathbb{E}^u \left[ V(\theta) + \int_{\mathfrak{s}}^{\theta} h(s, X, u_s) \, ds \Big| \mathcal{F}_{\mathfrak{s}} \right] \leq \operatorname*{essinf}_{v \in \mathfrak{U}} \mathbb{E}^v \left[ V(\mathfrak{t}) + \int_{\mathfrak{s}}^{\mathfrak{t}} h(s, X, v_s) \, ds \Big| \mathcal{F}_{\mathfrak{s}} \right],$$

the first in (56). [The reverse inequality holds on the event $\{\theta \leq \varrho_{\mathfrak{t}}\}$, as we know from (39).] The second inequality of (56) follows from the first, upon replacing $\theta$ by $\mathfrak{t}$, and $\mathfrak{t}$ by $\mathfrak{s}$.

Now (55) follows directly from (56), just by taking $\mathfrak{s} = 0$ there. $\square$

**7. A martingale characterization of saddle-points.** We are now in a position to provide necessary and sufficient conditions for the saddle-point property (10), in terms of appropriate martingales. These conditions are of obvious independent interest; they will also prove crucial when we try, in the next two sections, to prove constructively the existence of a saddle point $(u^*, \tau_*)$ for the stochastic game of control and stopping.

THEOREM 7.1.  *A pair $(u^*, \tau_*) \in \mathfrak{U} \times \mathcal{S}$ is a saddle point as in (10) for the stochastic game of control and stopping, if and only if the following three conditions hold:*

(i)  $g(X(\tau_*)) = V(\tau_*)$, *a.s.*
(ii)  $R^{u^*}(\cdot \wedge \tau_*)$ *is a $\mathbb{P}^{u^*}$-martingale; and*
(iii)  $R^u(\cdot \wedge \tau_*)$ *is a $\mathbb{P}^u$-submartingale, for every $u(\cdot) \in \mathfrak{U}$.*

The present section is devoted to the proof of this result. We shall derive first the conditions (i)–(iii) from the properties (10) of the saddle; then the reverse.

PROOF OF NECESSITY. Let us assume that the pair $(u^*, \tau_*) \in \mathfrak{U} \times \mathcal{S}$ is a saddle point for the game, that is, that the properties of (10) are satisfied.
• Using the definition of $\varrho_{\mathfrak{t}}$, the submartingale property $\mathbb{E}^{u^*}[R^{u^*}(\varrho_{\mathfrak{t}})|\mathcal{F}_{\mathfrak{t}}] \geq R^{u^*}(\mathfrak{t})$ from Proposition 5.3, the a.s. comparisons $Y^{u^*}(\tau_*) \leq R^{u^*}(\tau_*)$ and $Y^{u^*}(\varrho_{\tau_*}) = R^{u^*}(\varrho_{\tau_*})$, and the first property of the saddle in (10), we obtain

$$\mathbb{E}^{u^*}(Y^{u^*}(\tau_*)) \leq \mathbb{E}^{u^*}(R^{u^*}(\tau_*)) \leq \mathbb{E}^{u^*}(R^{u^*}(\varrho_{\tau_*}))$$
$$= \mathbb{E}^{u^*}(Y^{u^*}(\varrho_{\tau_*})) \leq \mathbb{E}^{u^*}(Y^{u^*}(\tau_*)).$$

But this gives, in particular, $\mathbb{E}^{u^*}(Y^{u^*}(\tau_*)) = \mathbb{E}^{u^*}(R^{u^*}(\tau_*))$, which, coupled with the earlier a.s. comparison, gives the stronger one $Y^{u^*}(\tau_*) = R^{u^*}(\tau_*)$, thus, also $g(X(\tau_*)) = V(\tau_*)$.
• Next, consider an *arbitrary* stopping rule $\tau \in \mathcal{S}$ with $0 \leq \tau \leq \tau_*$ and observe the string of inequalities

$$\mathbb{E}^{u^*}(R^{u^*}(\tau)) \leq \mathbb{E}^{u^*}(R^{u^*}(\varrho_\tau)) = \mathbb{E}^{u^*}(Y^{u^*}(\varrho_\tau))$$
$$\leq \mathbb{E}^{u^*}(Y^{u^*}(\tau_*)) = \mathbb{E}^{u^*}(R^{u^*}(\tau_*))$$

from Proposition 5.3, the definition of $\varrho_\tau$, the first property of the saddle, and property (i) just proved. On the other hand, from the second property of a saddle, from property (i) just proved and from the inequality (55), we get the second string of inequalities

$$\mathbb{E}^{u^*}(Y^{u^*}(\tau_*)) = \inf_{u \in \mathfrak{U}} \mathbb{E}^u(Y^u(\tau_*)) = \inf_{u \in \mathfrak{U}} \mathbb{E}^u(R^u(\tau_*))$$
$$\leq \inf_{u \in \mathfrak{U}} \mathbb{E}^u(R^u(\tau)) \leq \mathbb{E}^{u^*}(R^{u^*}(\tau)).$$

Combining the two strings, we deduce

$$(57) \quad \mathbb{E}^{u^*}(R^{u^*}(\tau)) = \inf_{u \in \mathfrak{U}} \mathbb{E}^u(R^u(\tau)) = \inf_{u \in \mathfrak{U}} \mathbb{E}^u(R^u(\tau_*)) = \mathbb{E}^{u^*}(R^{u^*}(\tau_*))$$

for *every* stopping rule $\tau \in \mathcal{S}$ with $0 \leq \tau \leq \tau_*$. This shows that $R^{u^*}(\cdot \wedge \tau_*)$ is a $\mathbb{P}^{u^*}$-martingale [cf. Exercise 1.3.26 in Karatzas and Shreve (1991)], and condition (ii) is established.
• It remains to show that, for *any* given $v(\cdot) \in \mathfrak{U}$, the process $R^v(\cdot \wedge \tau_*)$ is a $\mathbb{P}^v$-submartingale; equivalently, that for any stopping rules $\mathfrak{t}, \tau$ with $0 \leq \mathfrak{t} \leq \tau \leq \tau_*$, the inequality

$$(58) \quad \mathbb{E}^v\left[V(\tau) + \int_{\mathfrak{t}}^\tau h(s, X, v_s)\,ds \Big| \mathcal{F}_{\mathfrak{t}}\right] \geq V(\mathfrak{t}) \quad \text{holds a.s.}$$

Let us start by fixing a stopping rule $\tau$ as above, and recalling from (47) of Proposition 6.1 that

$$(59) \qquad \widehat{V}(\mathfrak{t};\tau) := \operatorname*{essinf}_{u \in \mathfrak{U}} \mathbb{E}^u \left[ V(\tau) + \int_{\mathfrak{t}}^{\tau} h(s, X, u_s) \, ds \Big| \mathcal{F}_{\mathfrak{t}} \right] \leq V(\mathfrak{t})$$

holds a.s. We'll be done, that is, we shall have proved (58), as soon as we have established that the reverse inequality

$$(60) \qquad \widehat{V}(\mathfrak{t};\tau) \geq V(\mathfrak{t}) \qquad \text{holds a.s.}$$

as well, for any given $\tau \in \mathcal{S}$ with $\mathfrak{t} \leq \tau \leq \tau_*$.

To this effect, let us consider for any $\varepsilon > 0$ the event $A_\varepsilon$ and the stopping rule $\theta_\varepsilon$ given as

$$A_\varepsilon := \{V(\mathfrak{t}) \geq \widehat{V}(\mathfrak{t};\tau) + \varepsilon\} \in \mathcal{F}_{\mathfrak{t}} \quad \text{and} \quad \theta_\varepsilon := \mathfrak{t} \cdot 1_{A_\varepsilon} + \tau \cdot 1_{A_\varepsilon^c},$$

respectively, and note $0 \leq \mathfrak{t} \leq \theta_\varepsilon \leq \tau \leq \tau_* \leq T$. From (57), we get

$$\mathbb{E}^{u^*}(R^{u^*}(\mathfrak{t})) = \mathbb{E}^{u^*}(R^{u^*}(\theta_\varepsilon)) = \mathbb{E}^{u^*}[R^{u^*}(\mathfrak{t}) \cdot 1_{A_\varepsilon} + R^{u^*}(\tau) \cdot 1_{A_\varepsilon^c}]$$

$$= \mathbb{E}^{u^*}[R^{u^*}(\mathfrak{t}) \cdot 1_{A_\varepsilon} + \mathbb{E}^{u^*}(R^{u^*}(\tau)|\mathcal{F}_{\mathfrak{t}}) \cdot 1_{A_\varepsilon^c}]$$

$$= \mathbb{E}^{u^*}\left[ V(\mathfrak{t}) \cdot 1_{A_\varepsilon} + \mathbb{E}^{u^*}\left(V(\tau) + \int_{\mathfrak{t}}^{\tau} h(s, X, u_s^*) \, ds \Big| \mathcal{F}_{\mathfrak{t}}\right) \cdot 1_{A_\varepsilon^c} \right.$$

$$\left. + \int_0^{\mathfrak{t}} h(s, X, u_s^*) \, ds \right]$$

$$\geq \mathbb{E}^{u^*}\left[ V(\mathfrak{t}) \cdot 1_{A_\varepsilon} + \widehat{V}(\mathfrak{t};\tau) \cdot 1_{A_\varepsilon^c} + \int_0^{\mathfrak{t}} h(s, X, u_s^*) \, ds \right]$$

$$\geq \varepsilon \cdot \mathbb{P}^{u^*}(A_\varepsilon) + \mathbb{E}^{u^*}\left[ \widehat{V}(\mathfrak{t};\tau) + \int_0^{\mathfrak{t}} h(s, X, u_s^*) \, ds \right].$$

That is,

$$(61) \qquad \mathbb{E}^{u^*}(R^{u^*}(\mathfrak{t})) - \varepsilon \cdot \mathbb{P}^{u^*}(A_\varepsilon) \geq \mathbb{E}^{u^*}\left[ \widehat{V}(\mathfrak{t};\tau) + \int_0^{\mathfrak{t}} h(s, X, u_s^*) \, ds \right].$$

As in (48), we write now the random variable $\widehat{V}(\mathfrak{t};\tau)$ of (59) in the form

$$\widehat{V}(\mathfrak{t};\tau) = \operatorname*{essinf}_{u \in \mathcal{U}_{[0,\mathfrak{t}]}^*} \mathbb{E}^u \left[ V(\tau) + \int_{\mathfrak{t}}^{\tau} h(s, X, u_s) \, ds \Big| \mathcal{F}_{\mathfrak{t}} \right]$$

$$= \lim_{k \to \infty} \mathbb{E}^{u^k}\left[ V(\tau) + \int_{\mathfrak{t}}^{\tau} h(s, X, u_s^k) \, ds \Big| \mathcal{F}_{\mathfrak{t}} \right]$$

for some sequence $\{u^k(\cdot)\}_{k \in \mathbb{N}}$ in $\mathcal{U}_{[0,\mathfrak{t}]}^*$, the set of admissible control strategies $u(\cdot) \in \mathfrak{U}$ that agree with $u^*(\cdot)$ a.e. on the stochastic interval $[[0, \mathfrak{t}]]$. Back into

(61), this gives

$$\mathbb{E}^{u^*}(R^{u^*}(\mathfrak{t})) - \varepsilon \cdot \mathbb{P}^{u^*}(A_\varepsilon)$$

$$\geq \mathbb{E}^{u^*}\left[\lim_k \mathbb{E}^{u^k}\left(V(\tau) + \int_{\mathfrak{t}}^{\tau} h(s, X, u_s^k)\,ds \Big| \mathcal{F}_{\mathfrak{t}}\right) + \int_0^{\mathfrak{t}} h(s, X, u_s^*)\,ds\right]$$

$$= \mathbb{E}^{u^*}\left[\lim_k \mathbb{E}^{u^k}\left(V(\tau) + \int_0^{\tau} h(s, X, u_s^k)\,ds \Big| \mathcal{F}_{\mathfrak{t}}\right)\right]$$

$$= \mathbb{E}^{u^*}\left[\lim_k \mathbb{E}^{u^k}(R^{u^k}(\tau)|\mathcal{F}_{\mathfrak{t}})\right]$$

$$= \lim_k \mathbb{E}^{u^*}[\mathbb{E}^{u^k}(R^{u^k}(\tau)|\mathcal{F}_{\mathfrak{t}})] \qquad \text{(bounded convergence)}$$

$$= \lim_k \mathbb{E}^{u^k}[\mathbb{E}^{u^k}(R^{u^k}(\tau)|\mathcal{F}_{\mathfrak{t}})] \qquad \text{[equation (24), Lemma 4.3]}$$

$$= \lim_k \mathbb{E}^{u^k}(R^{u^k}(\tau)) \geq \inf_{u \in \mathfrak{U}} \mathbb{E}^u(R^u(\tau)) = \mathbb{E}^{u^*}(R^{u^*}(\tau)) = \mathbb{E}^{u^*}(R^{u^*}(\mathfrak{t})).$$

The last claim follows from (57), the martingale property of (ii) that this implies, and $0 \leq \mathfrak{t} \leq \tau \leq \tau_*$. This shows $\mathbb{P}(A_\varepsilon) = 0$, and we get $V(t) < \widehat{V}(t; \tau) + \varepsilon$ a.s., for every $\varepsilon > 0$; letting $\varepsilon \downarrow 0$, we arrive at (60), and we are done. $\square$

PROOF OF SUFFICIENCY. Let us suppose now that the pair $(u^*, \tau_*) \in \mathfrak{U} \times \mathcal{S}$ satisfies the properties (i)–(iii) of Theorem 7.1; we shall deduce from them the properties of (10) for a saddle-point.

The $\mathbb{P}^u$-submartingale property of $R^u(\cdot \wedge \tau_*)$ in property (iii) gives $\mathbb{E}^u(R^u(\tau)) \leq \mathbb{E}^u(R^u(\tau_*))$ for all $u(\cdot) \in \mathfrak{U}$, thus, also

$$\inf_{u \in \mathfrak{U}} \mathbb{E}^u(R^u(\tau)) \leq \inf_{u \in \mathfrak{U}} \mathbb{E}^u(R^u(\tau_*)).$$

Taking here $\tau = 0$ and using the property (i) for $\tau_*$, as well as the $\mathbb{P}^{u^*}$-martingale property of $R^{u^*}(\cdot \wedge \tau_*)$ from (ii), we get

$$\inf_{u \in \mathfrak{U}} \mathbb{E}^u(Y^u(\tau_*)) = \inf_{u \in \mathfrak{U}} \mathbb{E}^u(R^u(\tau_*)) \geq R^u(0) = V = R^{u^*}(0)$$

$$= \mathbb{E}^{u^*}(R^{u^*}(\tau_*)) = \mathbb{E}^{u^*}(Y^{u^*}(\tau_*)).$$

Comparing the two extreme terms in this string, we obtain the second property of the saddle.

• We continue by considering stopping rules $\tau \in \mathcal{S}$ with $0 \leq \tau \leq \tau_*$. For such stopping rules, the fact that $R^{u^*}(\cdot \wedge \tau_*)$ is a $\mathbb{P}^{u^*}$-martingale [property (ii)] leads to

(62) $\quad Y^{u^*}(\tau) \leq R^{u^*}(\tau) = \mathbb{E}^{u^*}(R^{u^*}(\tau_*)|\mathcal{F}_\tau) = \mathbb{E}^{u^*}(Y^{u^*}(\tau_*)|\mathcal{F}_\tau) \qquad$ a.s.

and this gives the first property of the saddle for such stopping rules, upon taking expectations.

• Let us consider now stopping rules $\tau \in \mathcal{S}$ with $\tau_* \leq \tau \leq T$.
We shall establish for them the first property of the saddle, actually in the stronger form

$$\mathbb{E}^{u^*}(Y^{u^*}(\tau)|\mathcal{F}_{\tau_*}) \leq Y^{u^*}(\tau_*) \qquad \text{a.s.} \tag{63}$$

Now (63) is equivalent to

$$g(X(\tau_*)) \geq \mathbb{E}^{u^*}\left(g(X(\tau)) + \int_{\tau_*}^{\tau} h(t, X, u_t^*)\, dt \Big| \mathcal{F}_{\tau_*}\right)$$
$$= \mathbb{E}^{u^*}[Y^{u^*}(\tau_*, \tau)|\mathcal{F}_{\tau_*}],$$

a.s., for every $\tau \in \mathcal{S}_{\tau_*,T}$, thus to $g(X(\tau_*)) \geq Z^{u^*}(\tau_*)$, a.s. But from (19) and (21) the reverse of this inequality always holds, so (63) amounts to the requirement

$$g(X(\tau_*)) = Z^{u^*}(\tau_*) \qquad \text{a.s.} \tag{64}$$

To prove (64), recall from condition (ii) that $R^{u^*}(\cdot \wedge \tau_*)$ is a $\mathbb{P}^{u^*}$-martingale, and from (36) that it dominates $Y^{u^*}(\cdot \wedge \tau_*)$. But from Proposition 4.1, the process $Q^{u^*}(\cdot \wedge \tau_*)$ is the *smallest* $\mathbb{P}^{u^*}$-supermartingale that dominates $Y^{u^*}(\cdot \wedge \tau_*)$. Consequently, $R^{u^*}(\cdot \wedge \tau_*) \geq Q^{u^*}(\cdot \wedge \tau_*)$ and, equivalently, $V(\cdot \wedge \tau_*) \geq Z^{u^*}(\cdot \wedge \tau_*)$, hold a.s. But the reverse inequality also holds, thanks to the expression (26) for $V(\cdot)$, thus, in fact, $V(\cdot \wedge \tau_*) = Z^{u^*}(\cdot \wedge \tau_*)$, a.s. In particular, we get $V(\tau_*) = Z^{u^*}(\tau_*)$ a.s. Now (64) follows, in conjunction with condition (i).

• Finally, let us prove the first property of the saddle for an *arbitrary* stopping rule $\tau \in \mathcal{S}$. We start with the decomposition

$$\mathbb{E}^{u^*}(Y^{u^*}(\tau)) = \mathbb{E}^{u^*}(Y^{u^*}(\tau)1_{\{\tau \leq \tau_*\}} + Y^{u^*}(\tau)1_{\{\tau > \tau_*\}})$$
$$= \mathbb{E}^{u^*}(Y^{u^*}(\rho)1_{\{\tau \leq \tau_*\}} + Y^{u^*}(\nu)1_{\{\tau > \tau_*\}}),$$

where $\rho := \tau \wedge \tau_*$ belongs to $\mathcal{S}_{0,\tau_*}$ and $\nu := \tau \vee \tau_*$ is in $\mathcal{S}_{\tau_*,T}$. Thus, we have almost surely

$$Y^{u^*}(\rho) \leq \mathbb{E}^{u^*}(Y^{u^*}(\tau_*)|\mathcal{F}_\rho) \quad \text{and} \quad \mathbb{E}^{u^*}(Y^{u^*}(\nu)|\mathcal{F}_{\tau_*}) \leq Y^{u^*}(\tau_*),$$

from (62) and (63). Both events $\{\tau \leq \tau_*\}$, $\{\tau > \tau_*\}$ belong to $\mathcal{F}_\rho = \mathcal{F}_\tau \cap \mathcal{F}_{\tau_*}$, therefore,

$$\mathbb{E}^{u^*}(Y^{u^*}(\tau)) = \mathbb{E}^{u^*}(Y^{u^*}(\rho) \cdot 1_{\{\tau \leq \tau_*\}} + Y^{u^*}(\nu) \cdot 1_{\{\tau > \tau_*\}})$$
$$\leq \mathbb{E}^{u^*}(\mathbb{E}^{u^*}(Y^{u^*}(\tau_*)|\mathcal{F}_\rho) \cdot 1_{\{\tau \leq \tau_*\}} + \mathbb{E}^{u^*}(Y^{u^*}(\nu)|\mathcal{F}_\rho) \cdot 1_{\{\tau > \tau_*\}})$$
$$= \mathbb{E}^{u^*}(\mathbb{E}^{u^*}(Y^{u^*}(\tau_*) \cdot 1_{\{\tau \leq \tau_*\}}|\mathcal{F}_\rho) + \mathbb{E}^{u^*}(Y^{u^*}(\nu)|\mathcal{F}_{\tau_*}) \cdot 1_{\{\tau > \tau_*\}})$$
$$\leq \mathbb{E}^{u^*}(Y^{u^*}(\tau_*) \cdot 1_{\{\tau \leq \tau_*\}}) + \mathbb{E}^{u^*}(Y^{u^*}(\tau_*) \cdot 1_{\{\tau > \tau_*\}}) = \mathbb{E}^{u^*}(Y^{u^*}(\tau_*)).$$

This is the first property of the saddle in (10), established now for arbitrary $\tau \in \mathcal{S}$. □

**8. Optimality conditions for control.** We shall say that a given admissible control strategy $\widetilde{u}(\cdot) \in \mathfrak{U}$ is *optimal*, if it attains the infimum

$$(65) \qquad V = \inf_{v \in \mathfrak{U}} Z^v(0), \qquad \text{with } Z^v(0) = \sup_{\tau \in \mathcal{S}} \mathbb{E}^v[Y(\tau)].$$

Here and in what follows, we are using the notation of (19), (22) and (33). Clearly, if $(\widetilde{u}, \widetilde{\tau})$ is a saddle pair for the stochastic game, then $\widetilde{u}(\cdot)$ is an optimal control strategy.

THEOREM 8.1 (Necessary and sufficient condition for optimality of control). *A given admissible control strategy $u(\cdot) \in \mathfrak{U}$ is optimal, that is, attains the supremum in (65), if and only if it is thrifty, that is, satisfies*

$$(66) \qquad R^u(\cdot \wedge \tau_0^u) \qquad \text{is a } \mathbb{P}^u\text{-martingale}.$$

*And in this case, for every $0 \leq \varepsilon < 1$, we have in the notation of (33)*

$$(67) \qquad \tau_0^u(\varepsilon) = \varrho_0(\varepsilon) \qquad a.s.$$

PROOF OF SUFFICIENCY. Let us recall from (35) that $\tau_0^u(\varepsilon) \leq \tau_0^u$ holds a.s. for every $0 < \varepsilon < 1$, and from Proposition 4.2 that the process $Q^u(\cdot \wedge \tau_0^u)$ is a $\mathbb{P}^u$-martingale. Therefore, if $u(\cdot)$ is thrifty, we have

$$V \leq Z^u(0) = \mathbb{E}^u\left[Z^u(\tau_0^u(\varepsilon)) + \int_0^{\tau_0^u(\varepsilon)} h(s, X, u_s)\, ds\right]$$

$$\leq \mathbb{E}^u\left[\varepsilon + g(X(\tau_0^u(\varepsilon))) + \int_0^{\tau_0^u(\varepsilon)} h(s, X, u_s)\, ds\right]$$

$$\leq \varepsilon + \mathbb{E}^u\left[V(\tau_0^u(\varepsilon)) + \int_0^{\tau_0^u(\varepsilon)} h(s, X, u_s)\, ds\right]$$

$$= \varepsilon + \mathbb{E}^u[R^u(\tau_0^u(\varepsilon))] = \varepsilon + V.$$

In this string the second inequality comes from the definition of $\tau_0^u(\varepsilon)$ in (22); whereas the last equality is a consequence of thriftiness and of the inequality $\tau_0^u(\varepsilon) \leq \tau_0^u$. This gives the comparison $V \leq Z^u(0) \leq \varepsilon + V$ for every $0 < \varepsilon < 1$, therefore, $Z^u(0) = V$, the optimality of $u(\cdot)$. □

PROOF OF NECESSITY. Let us suppose now that $u(\cdot) \in \mathfrak{U}$ is optimal; we shall show that it is thrifty, and that (67) holds for every $0 \leq \varepsilon < 1$.
• We shall show first that, for this optimal $u(\cdot)$, we have $\tau_0^u = \varrho_0$ a.s., that is, (67) with $\varepsilon = 0$.

Let us observe that the $\mathbb{P}^u$-martingale property of $Q^u(\cdot \wedge \tau_0^u)$, coupled with the $\mathbb{P}^u$-submartingale property of $R^u(\cdot \wedge \varrho_0)$ from Proposition 5.3, and the a.s. inequality $\varrho_0 \leq \tau_0^u$ from (35) give

$$
\begin{aligned}
Z^u(0) - \mathbb{E}^u \int_0^{\varrho_0} h(s, X, u_s)\, ds \\
&= \mathbb{E}^u(Z^u(\varrho_0)) \\
&= \mathbb{E}^u[Z^u(\varrho_0) \cdot 1_{\{\tau_0^u = \varrho_0\}} + Z^u(\varrho_0) \cdot 1_{\{\tau_0^u > \varrho_0\}}] \\
&\geq \mathbb{E}^u[Z^u(\varrho_0) 1_{\{\tau_0^u = \varrho_0\}} + g(X(\varrho_0)) 1_{\{\tau_0^u > \varrho_0\}}] \\
&= \mathbb{E}^u[Z^u(\varrho_0) 1_{\{\tau_0^u = \varrho_0\}} + V(\varrho_0) 1_{\{\tau_0^u > \varrho_0\}}] \\
&\geq \mathbb{E}^u[V(\varrho_0) \cdot 1_{\{\tau_0^u = \varrho_0\}} + V(\varrho_0) \cdot 1_{\{\tau_0^u > \varrho_0\}}] \\
&= \mathbb{E}^u[V(\varrho_0)],
\end{aligned}
\tag{68}
$$

as well as

$$
\begin{aligned}
Z^u(0) &\geq \mathbb{E}^u\left[V(\varrho_0) + \int_0^{\varrho_0} h(s, X, u_s)\, ds\right] \\
&= \mathbb{E}^u[R^u(\varrho_0)] \geq R^u(0) = V.
\end{aligned}
\tag{69}
$$

We shall argue the validity of $\tau_0^u = \varrho_0$ by contradiction: we know from (35) that $\varrho_0 \leq \tau_0^u$ holds a.s., so let us assume

$$
\mathbb{P}^u(\tau_0^u > \varrho_0) > 0. \tag{70}
$$

Under the assumption (70), the first inequality in (68)—thus also in (69)—is *strict;* but this contradicts the optimality of $u(\cdot) \in \mathfrak{U}$. Thus, as claimed, we have $\tau_0^u = \varrho_0$ a.s. A similar argument leads to $\tau_0^u(\varepsilon) = \varrho_0(\varepsilon)$ a.s., for every $0 < \varepsilon < 1$, and (67) is proved.

• To see that this optimal $u(\cdot) \in \mathfrak{U}$ must also be *thrifty,* just observe that, as we have seen, equality prevails in (69); and that this, coupled with (67), gives $R^u(0) = \mathbb{E}^u[R^u(\tau_0^u)]$. It follows that the $\mathbb{P}^u$-submartingale $R^u(\cdot \wedge \varrho_0) \equiv R^u(\cdot \wedge \tau_0^u)$ is in fact a $\mathbb{P}^u$-martingale. □

The characterization of optimality presented in Theorem 8.1 is in the spirit of a similar characterization for optimal control with discretionary stopping in Dubins and Savage (1965) and in Maitra and Sudderth [(1996a), page 75]. In the context of these two sources, optimality amounts to the simultaneous validity of two conditions, "thriftiness" [i.e., condition (67)] and "equalization." In our context every control strategy is equalizing, so this latter condition becomes moot.

PROPOSITION 8.2. *If the admissible control strategy $u(\cdot) \in \mathfrak{U}$ is thrifty, then it is optimal; and the pair $(u, \tau_0^u) = (u, \varrho_0) \in \mathfrak{U} \times \mathcal{S}$ is then a saddle point for the stochastic game of control and stopping.*

PROOF. The first claim follows directly from Theorem 8.1. Now let us make a few observations:

(i) By the definition of $\varrho_0$ in (33) and the right-continuity of the process $V(\cdot)$, we have the a.s. equality $V(\varrho_0) = g(X(\varrho_0))$.

(ii) The process $R^u(\cdot \wedge \varrho_0)$ is a $\mathbb{P}^u$-martingale; this is because $u(\cdot)$, being optimal, must also be thrifty, as we saw in Theorem 8.1, and because $\varrho_0 = \tau_0^u$ holds a.s.

(iii) From Proposition 5.3, the process $R^v(\cdot \wedge \varrho_0)$ is a $\mathbb{P}^v$-submartingale, for every $v(\cdot) \in \mathfrak{U}$.

From these observations and Theorem 7.1, it is now clear that the pair $(u, \varrho_0)$ is a saddle point of the stochastic game. $\square$

**9. Constructing a thrifty control strategy and a saddle.** The theory of the previous section, culminating with Proposition 8.2, shows that in order to construct a saddle point for our stochastic game of control and stopping, all we need to do is find an admissible control strategy $u^*(\cdot) \in \mathfrak{U}$ which is *thrifty;* to wit, for which the condition (66) holds. Then the pair $(u^*, \tau_0^{u^*})$ will be a saddle point for our stochastic game.

To accomplish this, we shall start by assuming that, for each $(t, \omega)$, the mappings

(71) $\quad a \mapsto f(t, \omega, a) \quad \text{and} \quad a \mapsto h(t, \omega, a) \quad$ are continuous,

and that for the so-called *Hamiltonian* function

(72) $\quad H(t, \omega, a, p) := \langle p, \sigma^{-1}(t, \omega) f(t, \omega, a) \rangle + h(t, \omega, a),$

$t \in [0, T], \omega \in \Omega, a \in A, p \in \mathbb{R}^n$, the mapping $a \mapsto H(t, \omega, a, p)$ attains its infimum over the set $A$ at some $a^* \equiv \mathfrak{a}^*(t, \omega, p) \in A$, for any given $(t, \omega, p) \in [0, T] \times \Omega \times \mathbb{R}^n$, namely,

(73) $\quad \inf_{a \in A} H(t, \omega, a, p) = H(t, \omega, \mathfrak{a}^*(t, \omega, p), p).$

[This is the case, for instance, if the set $A$ is compact and the mapping $a \mapsto H(t, \omega, a, p)$ continuous.] Then it can be shown [see Lemma 1 in Beneš (1970), or Lemma 16.34 in Elliott (1982)] that the mapping $\mathfrak{a}^* : ([0, T] \times \Omega) \times \mathbb{R}^n \to A$ can be selected to be $(\mathcal{P} \otimes \mathcal{B}(\mathbb{R}^n)/\mathcal{A})$-measurable.

We shall deploy the martingale methodologies introduced in stochastic control in the seminal papers of Rishel (1970), Duncan and Varaiya (1971), Davis and Varaiya (1973) and Davis (1973), and presented in book form in Chapter 16 of Elliott (1982). The starting point of this approach is the observation that, for *every* admissible control strategy $u(\cdot) \in \mathfrak{U}$, the process

(74) $\quad R^u(\cdot \wedge \varrho_0) = V(\cdot \wedge \varrho_0) + \int_0^{\cdot \wedge \varrho_0} h(t, X, u_t) \, dt$

is a $\mathbb{P}^u$-submartingale with RCLL paths, and bounded uniformly on $[0,T] \times \Omega$; recall Propositions 5.3, 5.2 and Remark 5.5. This implies that the process $R^u(\cdot \wedge \varrho_0)$ admits a Doob–Meyer decomposition

$$(75) \qquad R^u(\cdot \wedge \varrho_0) = V + M^u(\cdot) + \Delta^u(\cdot).$$

Here $M^u(\cdot)$ is a uniformly integrable $\mathbb{P}^u$-martingale with RCLL paths and $M^u(0) = 0$, $M^u(\cdot) \equiv M^u(\varrho_0)$ on $[[\varrho_0, T]]$; the process $\Delta^u(\cdot)$ is predictable, with nondecreasing paths, $\Delta^u(T) \equiv \Delta^u(\varrho_0)$ integrable, and $\Delta^u(0) = 0$.

• A key observation now is that the $\mathbb{P}^u$-martingale $M^u(\cdot)$ can be represented as a stochastic integral, in the form

$$(76) \qquad M^u(\cdot) = \int_0^{\cdot} \langle \gamma(t), dW^u(t) \rangle.$$

Here $W^u(\cdot)$ is the $\mathbb{P}^u$-Brownian motion of (4), and $\gamma(\cdot)$ a predictable ($\mathcal{P}$-measurable) process that satisfies $\int_0^T \|\gamma(t)\|^2 \, dt < \infty$ and $\gamma(\cdot) \equiv 0$ on $[[\varrho_0, T]]$, a.s.

This is, of course, the predictable representation property of the filtration $\mathbb{F} = \{\mathcal{F}_t\}_{0 \le t \le T}$ [the augmentation of the filtration $\mathcal{F}_t^W = \sigma(W(s); 0 \le s \le t), 0 \le t \le T$, generated by the $\mathbb{P}$-Brownian motion $W(\cdot)$] under the equivalent change (5) of probability measure. For this result of Fujisaki et al. (1972), which is very useful in filtering theory, see Rogers and Williams (1987), pages 323 or Karatzas and Shreve (1998), Lemma 1.6.7. An important aspect of this representation is that *the same process $\gamma(\cdot)$ works for every $u(\cdot) \in \mathfrak{U}$ in* (76).

Next, let us take any two admissible control strategies $u(\cdot)$ and $v(\cdot)$ in $\mathfrak{U}$, and compare the resulting decompositions (75) on the stochastic interval $[[0, \varrho_0]]$. In conjunction with (74)–(76), (72) and (4), this gives

$$(77) \qquad \Delta^v(\cdot) - \Delta^u(\cdot) = \int_0^{\cdot} [H(t, X, v_t, \gamma(t)) - H(t, X, u_t, \gamma(t))] \, dt$$

on the interval $[[0, \varrho_0]]$. A brief, self-contained argument for the claims (76) and (77) is presented in the Appendix.

ANALYSIS. If we know that $\check{u}(\cdot) \in \mathfrak{U}$ is a thrifty control strategy, that is, the process $R^{\check{u}}(\cdot \wedge \tau_0^{\check{u}})$ is a $\mathbb{P}^{\check{u}}$-martingale, then $R^{\check{u}}(\cdot \wedge \varrho_0)$ is also a $\mathbb{P}^{\check{u}}$-martingale [just recall that we have $0 \le \varrho_0 \le \tau_0^{\check{u}}$ from (35)], thus $\Delta^{\check{u}}(\cdot) \equiv 0$ a.s. But then (77) gives

$$\Delta^v(\cdot) = \int_0^{\cdot} [H(t, X, v_t, \gamma(t)) - H(t, X, \check{u}_t, \gamma(t))] \, dt \qquad \text{on } [[0, \varrho_0]];$$

and because this process has to be nondecreasing for every admissible control strategy $v(\cdot) \in \mathfrak{U}$, we deduce the following necessary condition for thriftiness:

$$(78) \qquad H(t, X, \check{u}_t, \gamma(t)) = \inf_{a \in A} H(t, X, a, \gamma(t)) \qquad \text{a.e. on } [[0, \varrho_0]].$$

This is also known as the stochastic version of *Pontryagin's Maximum Principle*; cf. Kushner (1965), Haussmann (1986) and Peng (1990, 1993).

SYNTHESIS. The stochastic maximum principle of (78) suggests considering the admissible control strategy $u^*(\cdot) \in \mathfrak{U}$ defined by

$$(79) \qquad u_t^* = \begin{cases} \mathfrak{a}^*(t, X, \gamma(t)), & 0 \le t \le \varrho_0 \\ \mathrm{a}_\sharp, & \varrho_0 < t \le T \end{cases}$$

for an arbitrary but fixed element $\mathrm{a}_\sharp$ of the control set $A$. We are using here the "measurable selector" mapping $\mathfrak{a}^* : [0, T] \times \Omega \times \mathbb{R}^n \to A$ of (73).

With this choice, (77) leads to the comparison

$$\Delta^v(\cdot) = \Delta^{u^*}(\cdot) + \int_0^\cdot [H(t, X, v_t, \gamma(t)) - H(t, X, u_t^*, \gamma(t))] \, dt \ge \Delta^{u^*}(\cdot)$$

on the interval $[[0, \varrho_0]]$, therefore, also $R^v(\cdot \wedge \varrho_0) \ge V + M^v(\cdot) + \Delta^{u^*}(\cdot)$ from (75), for every $v(\cdot) \in \mathfrak{U}$. Taking expectations under $\mathbb{P}^v$, we obtain

$$0 \le \mathbb{E}^v[\Delta^{u^*}(\varrho_0)] \le \mathbb{E}^v[R^v(\varrho_0)] - V \qquad \forall \, v(\cdot) \in \mathfrak{U}.$$

But now we can take the infimum over $v(\cdot) \in \mathfrak{U}$ in the above string, and obtain

$$0 \le \inf_{v \in \mathfrak{U}} \mathbb{E}^v[\Delta^{u^*}(\varrho_0)] \le \inf_{v \in \mathfrak{U}} \mathbb{E}^v[R^v(\varrho_0)] - V = 0,$$

where the last equality comes from (55) and the sentence directly below it. We deduce

$$(80) \qquad \inf_{v \in \mathfrak{U}} \mathbb{E}^v[\Delta^{u^*}(\varrho_0)] = 0, \quad \text{thus also} \quad \Delta^{u^*}(\varrho_0) = 0 \qquad \text{a.s.}$$

from fairly standard weak compactness arguments, as in Davis (1973) page 592, Davis (1979) or Elliott (1982) pages 238–240.

• We follow now a reasoning similar to that used to prove (64) in Theorem 7.1: first, we note from (80) that

$$(81) \quad R^{u^*}(\cdot \wedge \varrho_0) = V + \int_0^\cdot \langle \gamma(t), dW^{u^*}(t) \rangle \qquad \text{is a } \mathbb{P}^{u^*}\text{-martingale},$$

and from (36) that it dominates $Y^{u^*}(\cdot \wedge \varrho_0)$. But from Proposition 4.1, the process $Q^{u^*}(\cdot \wedge \varrho_0)$ is the *smallest* $\mathbb{P}^{u^*}$-supermartingale that dominates $Y^{u^*}(\cdot \wedge \varrho_0)$. We deduce that $R^{u^*}(\cdot \wedge \varrho_0) \ge Q^{u^*}(\cdot \wedge \varrho_0)$ and, equivalently, $V(\cdot \wedge \varrho_0) \ge Z^{u^*}(\cdot \wedge \varrho_0)$, hold a.s. The reverse of this inequality also holds, thanks to the expression (26) for $V(\cdot)$, thus, in fact, $V(\cdot \wedge \varrho_0) = Z^{u^*}(\cdot \wedge \varrho_0)$.

In particular, we have almost surely, $Z^{u^*}(\varrho_0) = V(\varrho_0) = g(X(\varrho_0))$ (recall the definition of $\varrho_0$), thus, also $\tau_0^{u^*} \le \varrho_0$ from (22). Again, the reverse inequality holds, thanks now to (35), so, in fact, $\tau_0^{u^*} = \varrho_0$ holds a.s.

We conclude that the property (81) leads to the thriftiness condition (66) for the admissible control strategy $u^*(\cdot) \in \mathfrak{U}$ defined in (79). In conjunction with Proposition 8.2, this establishes the following existence and characterization result:

THEOREM 9.1. *Under the assumptions* (71)–(73) *of this section, the pair* $(u^*, \varrho_0) \in \mathfrak{U} \times \mathcal{S}$ *of* (79) *and* (33) *is a saddle point for the stochastic game, and we have* $\varrho_0 = \tau_0^{u^*}$ *a.s., in the notation of* (22). *Furthermore, the process* $V(\cdot \wedge \varrho_0)$ *is a continuous* $\mathbb{P}$-*semimartingale.*

Only the last claim needs discussion; from (74), (81) and (72), we get the representation

$$(82) \qquad V(t) = V - \int_0^t H(s, X, u_s^*, \gamma(s))\, ds + \int_0^t \langle \gamma(s), dW(s) \rangle$$

for $0 \le t \le \varrho_0$, and the claim follows.

This equation (82) can be written equivalently "backward," as

$$(83) \quad V(t) = g(X(\varrho_0)) + \int_t^{\varrho_0} H(s, X, u_s^*, \gamma(s))\, ds - \int_t^{\varrho_0} \langle \gamma(s), dW(s) \rangle$$

for $0 \le t \le \varrho_0$. Suitably modified to account for the constraint $V(\cdot) \ge g(X(\cdot))$, and with an appropriate definition for the "adjoint process" $\gamma(\cdot)$ on $[[\varrho_0, T]]$, the equation (83) can be extended to hold on $[[0, T]]$; this brings us into contact with the backward stochastic differential equation approach to stochastic games [Cvitanić and Karatzas (1996), Hamadène and Lepeltier (1995, 2000), Hamadène (2006)].

## APPENDIX

In order to make this paper as self-contained as possible, we shall present here a brief argument for the representation (76) of the $\mathbb{P}^u$-martingale $M^u(\cdot)$ in the Doob–Meyer decomposition (75), and for the associated identity (77).

We start with the "Bayes rule" computation

$$(84) \qquad \begin{aligned} M^u(t) &= \mathbb{E}^u[M^u(T)|\mathcal{F}_t] = \mathbb{E}^u[M^u(\varrho_0)|\mathcal{F}_t] \\ &= \frac{\mathbb{E}^u[\Lambda^u(\varrho_0) M^u(\varrho_0)|\mathcal{F}_t]}{\Lambda^u(t \wedge \varrho_0)} \end{aligned}$$

for $0 \le t \le T$ [e.g. Karatzas and Shreve (1991), page 193]; then the martingale representation property of the Brownian filtration (ibid., page 182) shows that the numerator of (84) can be expressed as the stochastic integral

$$(85) \quad N^u(t) := \mathbb{E}^u[\Lambda^u(\varrho_0) M^u(\varrho_0)|\mathcal{F}_t] = \int_0^t \langle \xi^u(s), dW(s) \rangle, \qquad 0 \le t \le T,$$

with respect to $W(\cdot)$, of some predictable process $\xi^u : [0,T] \times \Omega \to \mathbb{R}^n$ that satisfies $\xi^u(\cdot) \equiv 0$ a.e. on $[[\varrho_0, T]]$, and $\int_0^T \|\xi^u(t)\|^2\, dt < \infty$ a.s. We have recalled in (84) and (85) that $M^u(\cdot) \equiv M^u(\varrho_0)$ a.e. on $[[\varrho_0, T]]$, and $N^u(0) = M^u(0)\Lambda^u(0) = 0$.

On the other hand, for the exponential martingale of (3), we have the stochastic integral equation

$$(86) \qquad \Lambda^u(t \wedge \varrho_0) = 1 + \int_0^t \Lambda^u(s) \langle \varphi^u(s), dW(s) \rangle, \qquad 0 \leq t \leq T,$$

where we have set $\varphi^u(t) := \sigma^{-1}(t,X) f(t,X,u_t)$ for $0 \leq t \leq \varrho_0$, and $\varphi^u(t) := 0$ for $\varrho_0 < t \leq T$. Applying Itô's rule to the ratio $M^u(\cdot) = N^u(\cdot)/\Lambda^u(\cdot \wedge \varrho_0)$ of (84), in conjunction with (85), (86) and (4), we obtain then, for $0 \leq t \leq T$,

$$(87) \quad M^u(t) = \int_0^t \langle \gamma^u(s), dW^u(s) \rangle \qquad \text{where } \gamma^u(t) := \frac{\xi^u(t) - N^u(t)\varphi^u(t)}{\Lambda^u(t)}$$

is clearly predictable; it satisfies $\gamma^u(\cdot) \equiv 0$ a.e. on $[[\varrho_0, T]]$, as well as $\int_0^T \|\gamma^u(t)\|^2\, dt < \infty$ a.s.

• It remains to argue that the stochastic integrand of (87) does not depend on the admissible control process $u(\cdot) \in \mathfrak{U}$, as claimed in (76). Indeed, for arbitrary $u(\cdot) \in \mathfrak{U}$ and $u(\cdot) \in \mathfrak{U}$, we have

$$R^v(t \wedge \varrho_0) - \int_0^{t \wedge \varrho_0} [h(s,X,v_s) - h(s,X,u_s)]\, ds$$
$$= R^u(t \wedge \varrho_0) = V + \Delta^u(t) + \int_0^t \langle \gamma^u(s), dW^u(s) \rangle$$
$$= V + \Delta^u(t) + \int_0^t \langle \gamma^u(s), dW^v(s) \rangle$$
$$+ \int_0^t \langle \gamma^v(s), \varphi^v(s) \rangle\, ds - \int_0^t \langle \gamma^u(s), \varphi^u(s) \rangle\, ds, \qquad 0 \leq t \leq T.$$

Let us compare now this decomposition with the consequence

$$R^v(t \wedge \varrho_0) = V + \Delta^v(t) + \int_0^t \langle \gamma^v(s), dW^v(s) \rangle, \qquad 0 \leq t \leq T,$$

of (75) and (87). Identifying martingale terms, we see that $\gamma^u(\cdot) = \gamma^v(\cdot)$ holds a.e. on $[0,T] \times \Omega$, thus, (76) holds; identifying terms of bounded variation, we arrive at (77).

**Acknowledgments.** We are grateful to Professor William Sudderth for joint work and for his encouragement, both of which were critical in our undertaking this project. We are also indebted to the reviewers of the paper for their many incisive comments and for their corrections; they helped us improve this work very considerably.

Department of Mathematics  
Columbia University, MC 4438  
619 Mathematics Building  
New York, New York 10027  
USA  
E-mail: ik@math.columbia.edu

Department of Mathematics  
Baruch College, CUNY  
One Bernard Baruch Way, B6-230  
New York, New York 10010  
USA  
E-mail: ingrid-mona_zamfirescu@baruch.cuny.edu